\documentclass[12pt]{article}

\usepackage{amsmath,amssymb,amsthm}
\usepackage{epsfig}
\usepackage{arydshln}
\newcommand{\bfm}[1]{\mbox{\boldmath ${#1}$}}
\newcommand{\nonum}{\nonumber \\}
\newcommand\eq[1] {(\ref{#1})} 
\newcommand{\beqa}{\begin{eqnarray}}
\newcommand{\eeqa}[1]{\label{#1}\end{eqnarray}}
\newcommand{\beq}{\begin{equation}}
\newcommand{\eeq}[1]{\label{#1}\end{equation}}
\newcommand{\Grad}{\nabla}
\newcommand{\Div}{\nabla \cdot}

\newcommand{\Real}{\mathop{\rm Re}\nolimits}
    \newcommand{\Imag}{\mathop{\rm Im}\nolimits}
\newcommand{\Tr}{\mathop{\rm Tr}\nolimits}

\newcommand{\lang}{\langle}
\newcommand{\rang}{\rangle}
\newcommand{\Md}{\partial}

\newcommand{\Ga}{\alpha}
\newcommand{\Gb}{\beta}
\newcommand{\Gd}{\delta}
\newcommand{\Ge}{\epsilon}

\newcommand{\Gf}{\phi}

\newcommand{\Gg}{\gamma}

\newcommand{\Gk}{\kappa}

\newcommand{\Gl}{\lambda}
\newcommand{\Gn}{\eta}
\newcommand{\Gm}{\mu}

\newcommand{\Gt}{\theta}

\newcommand{\Gj}{\tau}

\newcommand{\GD}{\Delta}

\newcommand{\GO}{\Omega}

\newcommand{\BGe}{\bfm\epsilon}

\newcommand{\BGr}{\bfm\rho}

\newcommand{\BGs}{\bfm\sigma}

\newcommand{\BGj}{\bfm\tau}

\newcommand{\BGy}{\bfm\psi}

\newcommand{\BGD}{\bfm\Delta}

\newcommand{\BGS}{\bfm\Sigma}

\newcommand{\BGO}{\bfm\Omega}

\newcommand{\CT}{{\cal T}}

\newcommand{\bpm}{\begin{pmatrix}}
\newcommand{\epm}{\end{pmatrix}}
\newcommand\fig[1] {{\rm Figure}~\ref{fig:#1}}

\newcommand\labfig[1] {\label{fig:#1}}

\usepackage{url}
\usepackage{hyperref}

\def\Bc{{\bf c}}

\def\Be{{\bf e}}

\def\Bk{{\bf k}}

\def\Bn{{\bf n}}

\def\Bt{{\bf t}}
\def\Bu{{\bf u}}
\def\Bv{{\bf v}}
\def\Bw{{\bf w}}
\def\Bx{{\bf x}}

\def\BA{{\bf A}}
\def\BB{{\bf B}}
\def\BC{{\bf C}}

\def\BE{{\bf E}}

\def\BI{{\bf I}}

\def\BR{{\bf R}}
\def\BS{{\bf S}}
\def\BT{{\bf T}}

\usepackage{graphics}
\usepackage{subcaption}
\usepackage{amssymb}

\title{Planar polycrystals with extremal bulk and shear moduli}
\author{Graeme W. Milton}
\date{\small{Department of Mathematics, University of Utah, Salt Lake City, UT 84112, USA
\\Email: milton@math.utah.edu}}
\begin{document}
\maketitle
\vspace{2ex}
%******************************
\begin{abstract}
%******************************
  Here we consider the possible bulk and shear moduli of planar polycrystals built from a single crystal in various orientations.
  Previous work gave a complete characterization for crystals with orthotropic symmetry. Specifically, bounds were derived separately on the effective bulk and
  shear moduli, thus confining the effective moduli to lie within a rectangle in the (bulk, shear) plane. It was established that every point in this
  rectangle could be realized by an appropriate hierarchical laminate microgeometry, with the crystal taking different orientations in the layers, and the layers
  themselves being in different orientations.
  The bounds are easily extended to crystals with no special symmetry, but the path to constructing microgeometries that achieve every point in the
  rectangle defined by the bounds is considerably more difficult. We show that the two corners of the box having minimum bulk
  modulus are always attained by hierarchical laminates. For the other two corners we present
  algorithms for generating hierarchical laminates that attain them. Numerical evidence strongly suggests that the corner having maximum bulk and
  maximum shear modulus is always attained. For the remaining corner, with maximum bulk modulus and minimum shear modulus, 
   it is not yet clear whether the algorithm always succeeds, and hence
  whether all points in the rectangle are always attained. The microstructures we use are hierarchical laminate geometries
  that at their core have a self-similar microstructure, in the sense that the microstructure on one length scale is a rotation and rescaling of that on
  a smaller length scale.
\end{abstract}
\vspace{3ex}
%%%%%%%%%%%%%%%%%%%%%%%%%%%%%%%%%%%%%%%%%%%%%%%%%%%%%%%%%%%%%%%%%%%%%%%%
\section{Introduction}
\setcounter{equation}{0}
%%%%%%%%%%%%%%%%%%%%%%%%%%%%%%%%%%%%%%%%%%%%%%%%%%%%%%%%%%%%%%%%%%%%%%%%%%
This paper is a sequel to the work of Avellaneda et.al. \cite{Avellaneda:1996:CCP} where a complete characterization was given of
the possible bulk and shear moduli, $\Gk_*$ and $\Gm_*$, of planar polycrystals built from a single orthotropic crystal.
This was done by first finding bounds separately on $\Gk_*$ and $\Gm_*$ which thus define a rectangular box in the $(\Gk_*, \Gm_*)$-plane
of feasible moduli. Then polycrystal geometries were found that attain all points in the box. Our objective is to
extend this work to allow polycrystals built from crystals with elasticity tensors $\BC_0$
that are not orthotropic. The extension of the bounds is straightforward,
but the identification of the geometries that attain them is not. We will see that these geometries, and the associated proofs that they attain the bounds,
are highly non-trivial and quite different to those in \cite{Avellaneda:1996:CCP}. Moreover, unlike that in \cite{Avellaneda:1996:CCP}, the approach we take
completely avoids the difficulties of computing explicit expressions for the bounds and explicit expressions for the effective tensors of each microstructure
used in the construction. A side comment is that the geometries were first discovered by randomly constructing hierarchical laminate geometries and
finding those having effective tensors close to the bounds. 

The geometries are infinite rank laminates with a type of self-similar structure
as illustrated in \fig{1}. They resemble three-dimensional polycrystal geometries constructed by Nesi and Milton \cite{Nesi:1991:PCM} built
from a biaxial crystal and having the lowest possible isotropic effective conductivity (see also \cite{Tartar:1993:SRS,Bhattacharya:1994:RM} where related construction schemes occur in the context of convexification problems).
%The proofs encountered here require careful analysis and are considerably more difficult than those in either \cite{Avellaneda:1996:CCP} or \cite{Nesi:1991:PCM}.
A judicious choice of basis 
allows us to reformulate the problem as a problem of seeking trajectories in the complex plane, having
a prescribed form, that pass through a given point and which self-intersect (with the given point not being on the loop in the trajectory). From now on any use
of the word trajectory will imply that they have the prescribed form. 

\begin{figure}[ht!]
\includegraphics[width=0.9\textwidth]{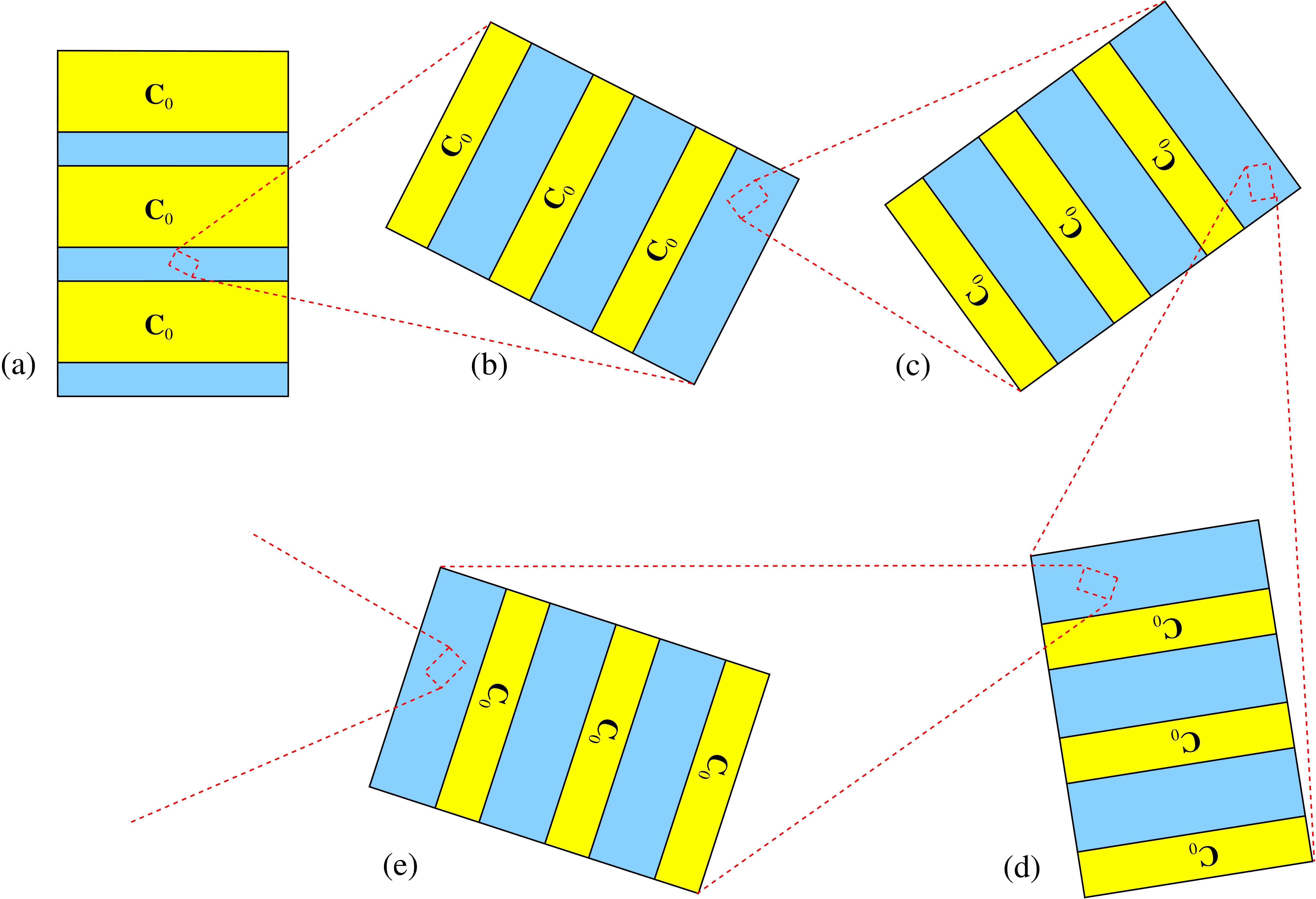}
\caption{A schematic illustration of the family of orthotropic microstructures attaining the bounds. By rotating the material and redefining $\BC_0$ to be a rotation of the original $\BC_0$, Here the blue material, that is not necessarily orthotropic, can be considered as a laminate
  of the original crystal and a rotation of the blue material, as in (b). Accordingly, we can iteratively replace the blue material by this laminate,
  as in (c), (d), (e), etc, and so the blue material itself can be regarding as a self-similar hierarchical polycrystal of the original crystal.
  The structure is self-similar in the sense that the structure at one length scale is a rescaling and rotation
  of the structure at other length scales. Ideally there should be a wide separation in length scales between the subsequent laminations. To obtain
  the desired isotropic material one takes the material in (a) and applies the construction scheme in \protect{\cite{Avellaneda:1996:CCP}}.
}
\labfig{1}
\end{figure}

The key idea
is to look for geometries such that the associated fields satisfy the attainability conditions for the bounds. This approach has been used time and again to construct
microgeometries attaining bounds derived from variational principles, or variational inequalities: for example, it was used in \cite{Milton:1981:CBT,Vigdergauz:1986:EEP,Milton:1988:VBE,Vigdergauz:1994:TDG,Grabovsky:1995:MMEa,Grabovsky:1995:MMEb,Gibiansky:2000:MCE,Sigmund:2000:NCE,Albin:2006:OCF,Albin:2007:MLE,Liu:2007:PIM} to find two-component and multicomponent microgeometries attaining the
Hashin-Shtrikman conductivity and bulk modulus
bounds; to show that ``coated laminate'' geometries with an isotropic effective elasticity tensor simultaneously attain the two-phase Hashin-Shtrikman bulk and shear modulus bounds \cite{Milton:1986:MPC} (as independently explicitly established in \cite{Norris:1985:DSE,Francfort:1986:HOB}); to show that
 infinite rank laminates attain the lower bound on the conductivity of an a three-dimensional polycrystal \cite{Nesi:1991:PCM}; to show that two-phase isotropic
geometries attaining the Hashin-Shtrikman conductivity bounds necessarily also attain the Hashin-Shtrikman bulk modulus bounds \cite{Grabovsky:1996:BEM}
(as also shown implicitly in \cite{Milton:1981:BEE} and explicitly in \cite{Berryman:1988:MRC})
and to extend this result to anisotropic composites of two possibly anisotropic phases; and to obtain laminate polycrystal geometries attaining the
Voigt and Reuss bounds on the effective bulk modulus \cite{Avellaneda:1989:OBE}. One popular approach to obtaining bounds from  variational principles is the translation method, or method of compensated compactness, of Tartar and Murat \cite{Tartar:1979:ECH,Murat:1985:CVH,Tartar:1985:EFC} and
Lurie and Cherkaev \cite{Lurie:1982:AEC,Lurie:1984:EEC} (see also the books \cite{Torquato:2001:RHM,Cherkaev:2000:VMS,Milton:2002:TOC,Tartar:2009:GTH}). Explicit attainability conditions on the fields for these bounds are presented in sections 25.3 and 25.4 of \cite{Milton:2002:TOC}. Following \cite{Avellaneda:1989:OBE} and \cite{Nesi:1991:PCM} we will
see that it suffices to look for fields meeting the attainability conditions and the appropriate differential constraints: the associated microgeometry then follows immediately
from the layout of these fields. 

Referring to \fig{2}, we will first prove attainability of points $A$ and $D$.  Due to the complexity of the problem, we do not yet have a proof that points $B$ and $C$ are always
attained, but we do have algorithms for finding geometries that attain them. Numerical evidence of Christian Kern presented here, based on the
algorithm for point $B$, strongly suggests that point $B$ is attained
for any choice of crystal elasticity tensor $\BC>0$. The question as to whether point $C$ is always attained using the corresponding algorithm has not been investigated, as this
requires exploring in a larger parameter space. In the case where all
four points $A, B, C,$ and $D$ are attained, the procedure for finding microgeometries that attain any point within the rectangle is exactly the same as that given in Section 4.1 of
\cite{Avellaneda:1996:CCP}. These arguments additionally imply that if only the points $A$, $D$, and $B$ are attained, then all points along the line segment joining $A$ and $D$
and  along the line segment joining $A$ and $B$ are also attained. 

\begin{figure}[ht!]
\includegraphics[width=0.9\textwidth]{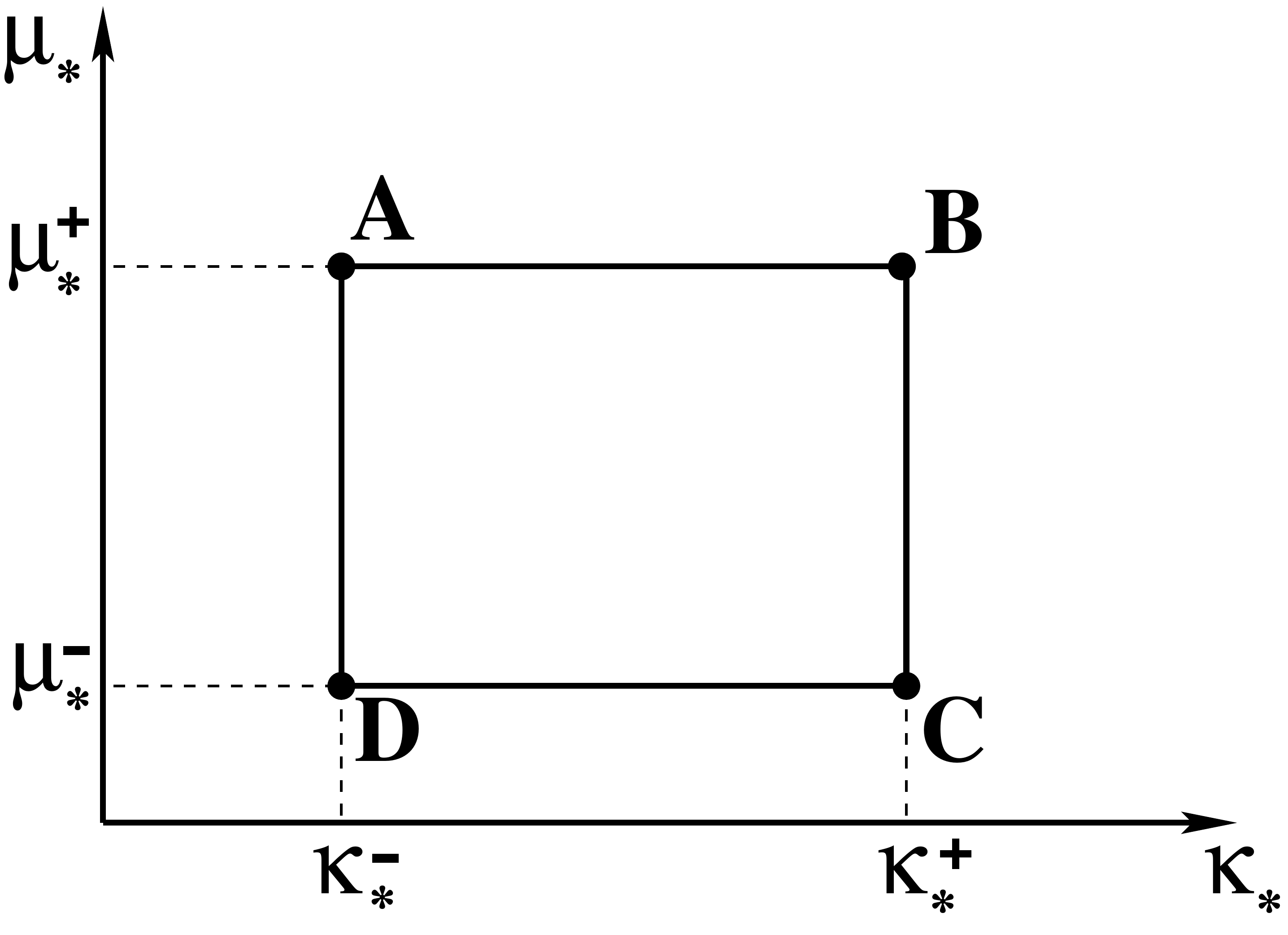}
\caption{The bounds on the possible pairs $(\Gk_*,\Gm_*)$ of the bulk modulus $\Gk_*$ and shear modulus $\Gm_*$ of isotropic planar polycrystals are given by
  the rectangular box, $\Gk_*^-\leq\Gk_*\leq\Gk_*^+, \Gm_*^-\leq\Gm_*\leq\Gm_*^+$. We show that the points $A$ and $D$ are attained, and we provide algorithms for
  finding microstructures that attain the points $B$ and $C$, though it is not yet clear that these algorithms always work.}
\labfig{2}
  \end{figure}

Part of the motivation for considering this problem is the renewed interest \cite{Huang:2011:TDM,Andreassen:2014:DME,Xia:2015:DMT,Berger:2017:MMT,Ostanin:2017:PCC,Yera:2020:TDE} in the range of possible bulk and shear moduli of isotropic composites of
two isotropic phases, both for planar elasticity and in the full three dimensional case, adding to earlier works reviewed in the books
\cite{Cherkaev:2000:VMS,Milton:2002:TOC,Allaire:2002:SOH,Torquato:2002:RHM,Tartar:2009:GTH}.
(There has also been work on the range of effective elasticity tensors for anisotropic composites, as surveyed in these books, with recent advances in the papers
\cite{Milton:2016:PEE,Milton:2016:TCC,Milton:2018:NOP}).
This interest has been driven by the advent of 3d-printing that allows one
to produce tailor made microstructures. In exploring this range it is easier to produce structures with anisotropic effective elasticity tensors. One can
convert them into microstructures with isotropic elasticity tensors by making polycrystals of the anisotropic structures (with the length scales in the
polycrystal being much larger than the length scales in the anisotropic structures so that one can apply homogenization theory). In doing this one
is faced with the question of determining the possible effective bulk and shear moduli of the polycrystal.

%%%%%%%%%%%%%%%%%%%%%%%%%%%%%%%%%%%%%%%%%%%%%%%%%%%%%%%%%%%%%%%%%%%%%%%%
\section{Preliminaries}
\setcounter{equation}{0}
%%%%%%%%%%%%%%%%%%%%%%%%%%%%%%%%%%%%%%%%%%%%%%%%%%%%%%%%%%%%%%%%%%%%%%%%%%
The planar linear elasticity equations in a periodic polycrystal take the form
\beq \BGs(\Bx)=\BC(\Bx)\BGe(\Bx),\quad\Div\BGs=0,\quad \BGe=[\Grad\Bu+(\Grad\Bu)^T]/2,
\eeq{1.a}
where $\BGs(\Bx), \BGe(\Bx)$ and $\Bu(\Bx)$ are the stress, strain, and displacement field,
and the elasticity fourth order tensor field $\BC(\Bx)$ with elements in Cartesian coordinates given by
\beq C_{ijk\ell}(\Bx)=R_{ia}(\Bx)R_{jb}(\Bx)R_{kc}(\Bx)R_{\ell d}\{\BC_0\}_{abcd}, \eeq{1.b}
where $\BC_0$ is the elasticity tensor of the original crystal, with elements $\{\BC_0\}_{abcd}$  and the periodic rotation field
$\BR(\Bx)$, with elements $R_{\Ga\Gb}(\Bx)$  determines its orientation throughout the polycrystal. If we look for solutions
where $\BGs(\Bx)$ and $\BGe(\Bx)$ have the same periodicity as the composite, which constitutes the ``cell-problem''
in periodic homogenization, then their
average values are linearly related, and this linear relation:
\beq \lang\BGs\rang=\BC_*\lang\BGe\rang, \eeq{1.c}
determines the effective elasticity tensor $\BC_*$. Here, and thereafter, the angular brackets $\lang\cdot\rang$ denote a volume average
over the unit cell of periodicity. If we test the material against two applied
shears
\beq \lang\BGe^{(1)}\rang=\frac{1}{2}\bpm 1 & 0 \\ 0 & -1 \epm, \quad \lang\BGe^{(2)}\rang=\frac{1}{2}\bpm 0 & 1 \\ 1 & 0 \epm,
\eeq{1.d}
then there will be two resulting stress, strain, and displacement fields, $\BGs^{(i)}(\Bx), \BGe^{(i)}(\Bx)$ and $\Bu^{(i)}(\Bx)$, $i=1,2$. 
This motivates the introduction of complex fields
\beq \BGs(\Bx)=\BGs^{(1)}(\Bx)+i\BGs^{(2)}(\Bx), \quad \BGe(\Bx)=\BGe^{(1)}(\Bx)+i\BGe^{(2)}(\Bx), \quad \Bu(\Bx)=\Bu^{(1)}(\Bx)+i\Bu^{(2)}(\Bx),
\eeq{1.e}
and then the equations \eq{1.a} will still hold with the elasticity tensor being real.

Using the elasticity equations and integration by parts, one sees that
\beq \overline{\lang\BGe\rang}:\BC_*\lang\BGe\rang=\lang\overline{\BGe}:\BC\BGe\rang,
\quad 
\overline{\lang\BGs\rang}:\BS_*\lang\BGs\rang=\lang\overline{\BGs}:\BS\BGs\rang,
\eeq{1.f}
where the overline denotes complex conjugation, $\BS(\Bx)=[\BC(\Bx)]^{-1}$ and $\BS_*=[\BC_*]^{-1}$ are the local and effective compliance tensors
(and the inverse is on the space of symmetric matrices) and the colon $:$ denotes a double contraction of indices,
\beq \BA :\BB\equiv\sum_{i=1,2}\sum_{j=1,2}A_{ij}B_{ij}.
\eeq{1.g}

We take an orthonormal basis
\beq \BGr_1=\frac{1}{2}\bpm 1 & i \\ i & -1 \epm, \quad
 \BGr_2=\frac{1}{2}\bpm 1 & -i \\ -i & -1 \epm, \quad 
 \BGr_3=\frac{1}{2}\bpm 1 & i \\ -i & 1 \epm, \quad
 \BGr_4=\frac{1}{2}\bpm 1 & -i \\ i & 1 \epm.
 \eeq{1.1}
 In this basis, the identity matrix $\BI=\BGr_3+\BGr_4$ is represented by the vector
 \beq \BI=\bpm 0 \\ 0 \\ 1 \\ 1 \epm. \eeq{1.1a}
 Under a rotation $\BR_\Gt$ anticlockwise by an angle $\Gt$ the matrix $\BGr_1$ transforms to
 \beq \BGr_1'=\BR_\Gt\BGr_1\BR_\Gt^T=\frac{1}{2}\bpm \cos\Gt & -\sin\Gt \\ \sin\Gt & \cos\Gt \epm\bpm 1 & i \\ i & -1 \epm\bpm \cos\Gt & \sin\Gt \\ -\sin\Gt & \cos\Gt \epm
 =e^{i2\Gt}\BGr_1, \eeq{1.2}
 and by taking complex conjugates we see that $\BGr_2$ transforms to $\BGr_2'=e^{-i2\Gt}\BGr_2$, while $\BGr_3$, and $\BGr_4$ are clearly rotationally invariant. Thus, in the basis
 \eq{1.1}, $\BR_\Gt$ is represented by the matrix
 \beq \BR_\Gt=\bpm  e^{i2\Gt} & 0 & 0 & 0 \\ 0 & e^{-i2\Gt} & 0 & 0 \\ 0 & 0 & 1 & 0\\ 0 & 0 &0 & 1 \epm.
 \eeq{1.3}
 An isotropic effective elastic tensor is represented by the matrix
 \beq  \BC_*=\bpm  2\mu_* & 0 & 0 & 0 \\ 0 & 2\mu_* & 0 & 0 \\ 0 & 0 & \Gk_* & \Gk_* \\ 0 & 0 & \Gk_* & \Gk_* \epm.
 \eeq{1.10}
 Note that the tensors $\BC_0$ and $\BS_0=\BC_0^{-1}$ are represented by Hermitian rather than real matrices, because the basis \eq{1.1} consists of
 complex matrices.

 If we layer two materials in direction $\Bn=(0,1)$ then the differences $\BGs_1-\BGs_2$ in stress, $\BE_1-\BE_2$ in displacement gradient, and $\BGe_1-\BGe_2$
 in strain across the interface must satisfy the jump conditions:
 \beq (\BGs_1-\BGs_2)\Bn=0,\quad (\BE_1-\BE_2)^T\Bt=0,\quad \Bt\cdot(\BGe_1-\BGe_2)\Bt=0,
 \eeq{1.11}
 where $\Bt=(1,0)$ is parallel to the layer interface. The transpose in \eq{1.11} arises because we choose the notation (contrary to that commonly used in continuum
 mechanics) where  $\BE=\Grad\Bu$ has elements $E_{ij}=\Md_iu_j$.

 This implies, that in the basis \eq{1.1}, these jumps must be of the form
 \beqa \BGs_1-\BGs_2 & = & \Ga_1\bpm 1 \\ 1 \\ 1 \\ 1 \epm, \quad \BE_1-\BE_2=\bpm \Gb_1 \\ \Gb_2 \\ -\Gb_1 \\ -\Gb_2 \epm, \nonum
 \BGe_1-\BGe_2 & = & \bpm \Gg_1 \\ \Gg_2 \\ \Gg_3 \\ \Gg_3 \epm\quad\text{with}\quad \Gg_1+ \Gg_2+ 2\Gg_3=0.
 \eeqa{1.12}
 Now consider a periodic field $\BE(\Bx)=\Grad\Bu(\Bx)$. Letting $\underline{\Bu}$ denote the periodic field $\Bu-\Bx^T\lang\BE\rang$,
 for $\Bk\ne 0$, $\BE$ has Fourier components
 \beqa \widehat{\BE}(\Bk) & = & i\Bk\otimes\underline{\Bu}=\bpm ik_1\widehat{\underline{u}}_1(\Bk) & ik_1\widehat{\underline{u}}_2(\Bk) \\ ik_2\widehat{\underline{u}}_1(\Bk) & ik_2\widehat{\underline{u}}_2(\Bk) \epm \nonum
 & = & \Ga[k_1(\BGr_1+\BGr_3)-ik_2(\BGr_1-\BGr_3)]+\Gb[k_1(\BGr_2+\BGr_4)+ik_2(\BGr_2-\BGr_4)], \nonum
 & ~ &
 \eeqa{1.13}
 where
 \beq \Ga= (i\widehat{\underline{u}}_1(\Bk)+\widehat{\underline{u}}_2(\Bk))/2,\quad \Gb= (i\widehat{\underline{u}}_1(\Bk)-\widehat{\underline{u}}_2(\Bk))/2.
\eeq{1.13a}
 So, for $\Bk\ne 0$, $\widehat{\BE}(\Bk)$ is represented by the vector
 \beq \widehat{\BE}(\Bk) = \bpm \Ga\overline{\Bk} \\ \Gb\Bk \\ \Ga\Bk \\ \Gb\overline{\Bk} \epm,
 \eeq{1.14}
 in which $\Bk=k_1+ik_2$ and $\overline{\Bk}=k_1-ik_2$. Here and thereafter the overline will denote complex conjugation.

 On the other hand with $\BR_\perp$ denoting the matrix for a $90^\circ$ rotation, if a periodic field
 $\BGj(\Bx)$ is divergence free, $\Div\BGj=0$  then $\BR_\perp\BGj$ has zero curl, and so $\BR_\perp(\BGj-\lang\BGj\rang)$ can be expressed as the gradient of a
 periodic potential that we label as $\BGy$. Thus, for $\Bk\ne 0$, the Fourier components take the form
 \beqa \widehat{\BGj}(\Bk) & = & \bpm k_2\widehat{\psi}_1(\Bk) & k_2\widehat{\psi}_2(\Bk) \\ -k_1\widehat{\psi}_1(\Bk) & -k_1\widehat{\psi}_2(\Bk) \epm \nonum
 & = & \Gg[k_2(\BGr_1+\BGr_3)+ik_1(\BGr_1-\BGr_3)]+\Gd[k_2(\BGr_2+\BGr_4)-ik_1(\BGr_2-\BGr_4)], \nonum
 &~&
 \eeqa{1.15}
 where
 \beq \Gg=(\widehat{\psi}_1(\Bk)-i\widehat{\psi}_2(\Bk))/2,\quad \Gd=(\widehat{\psi}_1(\Bk)+i\widehat{\psi}_2(\Bk))/2.
\eeq{1.15a}
The Fourier component $\widehat{\BGj}(0)$ can be arbitrary. So, for $\Bk\ne 0$, $\widehat{\BGj}(\Bk)$ is represented by the vector
 \beq \widehat{\BGj}(\Bk) = i\bpm \Gg\overline{\Bk} \\ -\Gd\Bk \\ -\Gg\Bk \\ \Gd\overline{\Bk} \epm.
 \eeq{1.16}
 %%%%%%%%%%%%%%%%%%%%%%%%%%%%%%%%%%%%%%%%%%%%%%%%%%%%%%%%%%%%%%%%%%%%%%%%
\section{Bounds and their attainability}
\setcounter{equation}{0}
%%%%%%%%%%%%%%%%%%%%%%%%%%%%%%%%%%%%%%%%%%%%%%%%%%%%%%%%%%%%%%%%%%%%%%%%%%
The lower bound on the bulk modulus is the Reuss-Hill bound \cite{Hill:1952:EBC},
\beq 1/\Gk_*\leq k\equiv\Tr(\Be)\quad \text{where}\quad \Be=\BS_0\BI,
\eeq{2.30}
and $\Be$ is symmetric because the range of $\BS_0$ consists of symmetric matrices. 
Defining $\Gk_*^-=1/k$ we then have the lower bound $\Gk_*\geq\Gk_*^-$.
More generally, allowing for anisotropic composites,
we have the bound $\BI:\BS_*\BI\leq k$. This bound will be attained \cite{Avellaneda:1989:OBE} if we can find a rotation field $\BR_{\Gt(\Bx)}$ such that
\beq \BGe(\Bx)=\BR_{\Gt(\Bx)}\Be \eeq{2.30a}
is a periodic strain field, i.e., the symmetrized gradient of a vector field.
Then we take as our polycrystal the material with in the basis \eq{1.1} an elasticity tensor
\beq \BC(\Bx)=\BR_{\Gt(\Bx)}\BC_0[\BR_{\Gt(\Bx)}]^\dagger, \eeq{2.20}
where the dagger denotes the complex conjugate of the transpose. This ensures that
\beq \BC(\Bx)\BGe(\Bx)=\BR_{\Gt(\Bx)}\BC_0[\BR_{\Gt(\Bx)}]^\dagger\BR_{\Gt(\Bx)}\Be=\BR_{\Gt(\Bx)}\BC_0\Be=\BR_{\Gt(\Bx)}\BI=\BI.
\eeq{2.30b}
In other words, the elasticity equations are solved with a stress $\BI$ that is constant.

Now suppose we can find a fourth-order tensor $\BT$, called the translation, such that $\BT$ is self-adjoint and
\beq \lang\overline{\BE}:\BT\BE\rang\geq \overline{\lang\BE\rang}:\BT\lang\BE\rang \eeq{2.30c}
for all periodic displacement gradients $\BE=\Grad\Bu$. Then if $\BC(\Bx)\geq\BT$ holds on the space of all complex matrices and for all $\Bx$,
the first identity in \eq{1.f} implies
\beq
\overline{\lang\BE\rang}:\BC_*{\lang\BE\rang}=\lang\overline{\BE}:\BC\BE\rang\geq\lang\overline{\BE}:\BT\BE\rang\geq\overline{\lang\BE\rang}:\BT\lang\BE\rang,
\eeq{2.30d}
where $\BE=\Grad\Bu$ is periodic, and $\Bu(\Bx)$ solves the elasticity equations \eq{1.a}. 
Note that in \eq{1.f} we have replaced $\BGe$ by $\BE=\Grad\Bu$ since $\BC$ and $\BC_*$ annihilate the antisymmetric parts of
$\BE$ and $\lang\BE\rang$. Thus we are left with the inequality $\BC_*-\BT\geq 0$, which is the comparison bound \cite{Avellaneda:1988:ECP}, a special
case of the translation method, or method of compensated compactness, of Tartar and Murat \cite{Tartar:1979:ECH,Murat:1985:CVH,Tartar:1985:EFC} and
Lurie and Cherkaev \cite{Lurie:1982:AEC,Lurie:1984:EEC} for bounding the effective tensors of composites.

Similarly, suppose one can find a translation $\BT$ such that $\BT$ is self-adjoint and
\beq \lang\overline{\BGs}:\BT\BGs\rang\geq \overline{\lang\BGs\rang}:\BT\lang\BGs\rang, \eeq{2.30e}
for all periodic fields $\BGs$ such that $\Div\BGs=0$. Then if, for all $\Bx$, $\BS(\Bx)\geq\BT$ on the space of symmetric matrices,
the second identity in \eq{1.f} implies
\beq
\overline{\lang\BGs\rang}:\BS_*\lang\BGs\rang=\lang\overline{\BGs}:\BS\BGs\rang\geq\lang\overline{\BGs}:\BT\BGs\rang\geq\overline{\lang\BGs\rang}:\BT\lang\BGs\rang,
\eeq{2.30f}
where the periodic stress field $\BGs(\Bx)$ solves the elasticity equations \eq{1.a}. Thus on the space of symmetric matrices
we are left with the comparison bound $\BS_*-\BT\geq 0$.

Observe that \eq{2.30c} is satisfied as an equality if $\BT$ maps displacement gradients to divergence free fields, and
\eq{2.30e} is satisfied as an equality if, conversely, $\BT$ maps  divergence free fields to displacement gradients. In
these cases the quadratic form associated with $\BT$ is a null-Lagrangian \cite{Murat:1978:CPC,Ball:1981:NLW}.

An upper bound on the effective bulk modulus is obtained using the translation represented by the matrix
\beq  \BT_0=\bpm -t_0 & 0 & 0 & 0 \\ 0 & -t_0 & 0 & 0 \\ 0 & 0 & t_0 & 0\\ 0 & 0 & 0 & t_0 \epm,
\eeq{2.31}
where $t_0$ is real. This is rotationally invariant. A side remark is that $\BA:\BT\BA$ (not $\overline{\BA}:\BT\BA$ which is real for real $t_0$) equals
$2t_0\det[\BR_\perp^T\BA\BR_\perp$ for any complex matrix $\BA$.
One can easily check that $\BT_0$ maps fields of the form \eq{1.16} to those of the form \eq{1.14} and thus maps
divergence free fields to displacement gradients, ensuring equality in \eq{2.30e}. We choose $t_0>0$ so that
\beq \BS_0-\BT_0\geq 0,\quad \det(\BI-\BC_0\BT_0)=0, \eeq{2.32}
and then the upper bound on the effective bulk modulus is obtained from the inequality $\BS_*-\BT_0\geq 0$ which implies
\beq \Gk_*\leq \Gk_*^+\equiv 1/(2t_0). \eeq{2.33}
In practice we accomplish this taking $t_0$ to be the lowest positive root of the equation $\det(\BI-\BC_0\BT_0)=0$.
There then exists a vector $\Bc$ such that
\beq \Bc=\BC_0\BT_0\Bc, \eeq{2.34a}
in which $\Bc$ and $\BT_0\Bc$ represent real symmetric matrices because $\Bc$ is in the range of $\BC_0$ and can be taken
to be real because $\BT_0$ represents a real fourth order tensor (thus, if $\Bc$ satisfies \eq{2.34a} so does
$\overline{\Bc}$ and also $\Bc+\overline{\Bc}$, the latter being real). 
Then $c_2=\overline{c_1}$ and $c_3=c_4$ is real, and if the latter is nonzero we may normalize $\Bc$ so
that $c_3=c_4=1$. Then the inequality that
\beq 0\leq\overline{\Bc}\cdot\BS_0\Bc=\overline{\Bc}\cdot\BT_0\Bc=2t_0(1-|c_1|^2) \eeq{2.35}
implies $|c_1|\leq 1$.
The bound $\BS_*-\BT\geq 0$ will be attained if we can find a divergence free stress field $\BGs(\Bx)$ that in the basis \eq{1.1} takes the form
\beq \BGs(\Bx)=\Gg(\Bx)\BR_{\Gt(\Bx)}\Bc, \eeq{2.34}
and is such that $\lang\BGs\rang\ne 0$.
Then by the properties of $\BT$, $\BE(\Bx)=\BT\BGs(\Bx)$ is a displacement gradient and 
\beq \BE(\Bx)=\BT\BGs(\Bx)=\Gg(\Bx)\BR_{\Gt(\Bx)}\BT\Bc. \eeq{2.8}
Taking a polycrystal with in the basis \eq{1.1} an elasticity tensor
\beq \BC(\Bx)=\BR_{\Gt(\Bx)}\BC_0[\BR_{\Gt(\Bx)}]^\dagger, \eeq{2.9}
we see, using \eq{2.34a}, that
\beq \BC(\Bx)\BE(\Bx)=\Ga(\Bx)\BC(\Bx)\BR_{\Gt(\Bx)}\BT\Bc=\Gg(\Bx)\BR_{\Gt(\Bx)}\BC_0\BT\Bc=\Gg(\Bx)\BR_{\Gt(\Bx)}\Bc=\BGs(\Bx).
\eeq{2.10}
As the fields $\BE(\Bx)$ and $\BGs(\Bx)$ solve the constitutive relation, and the differential constraints, their averages are
related by the effective elasticity tensor $\BC_*$:
\beq \lang\BGs\rang=\BC_*\lang\BE\rang, \eeq{2.11}
and taking averages of both sides of $\BE(\Bx)=\BT\BGs(\Bx)$ gives
\beq \lang\BE\rang=\BT\lang\BGs\rang=\BT\BC_*\lang\BE\rang. \eeq{2.12}
So the bound $\BS_*-\BT\geq 0$ is attained.

Bounds on the effective shear modulus are obtained using the translation represented by the matrix
\beq \BT=\bpm t_1 & 0 & 0 & 0 \\ 0 & t_2 & 0 & 0 \\ 0 & 0 & -t_1 & 0\\ 0 & 0 & 0 & -t_2 \epm,
\eeq{2.1}
which is rotationally invariant. This translation maps fields of the form \eq{1.14} to those of the form \eq{1.16} and vice-versa,
so both \eq{2.30c} and \eq{2.30e} hold as equalities.

First, to obtain an upper bound on $\Gm_*$
we require that $t_1$ and $t_2$ be chosen such that
\beq \BS_0-\BT\geq 0,\quad \det(\BI-\BC_0\BT)=0. \eeq{2.2}
Among pairs $(t_1,t_2)$ satisfying this, we select the pair having the maximal value of $t_1$, since the bounds $\BS_*-\BT\geq 0$ imply
\beq 1/(2\Gm_*)\geq t_1,\quad 1/(2\Gm_*)\geq t_2, \quad 1/(2\Gm_*)\geq t_1,\quad 1/(2\Gm_*)\geq t_2,\quad t_1+t_2+1/\Gk_*\geq 0. \eeq{2.2a}
Defining this maximum value of $t_1$ to be $1/(2\Gm_*^+)$ we then have the upper bound $\Gm_*\leq\Gm_*^+$.
In practice we find $\Gm_*^+$  by taking in the $(t_1,t_2)$ plane the simply connected loop of the curve $\det(\BI-\BC_0\BT)=0$ that encloses the origin,
and surrounds it, and finding the maximal value of $t_1$ along this curve. To justify this procedure observe that 
\beq \CT=\{(t_1,t_2)\in\mathbb{R}^2 | \BS_0-\BT(t_1,t_2)\geq 0 \quad \text{on symmetric matrices}\}
\eeq{2.2aaaa}
is a convex compact set. The origin $(0,0)$ is obviously in the interior of $\CT$. Thus, $\Md\CT$ is the boundary of an open and bounded convex set
containing the origin.
Every point on $\Md\CT$ satisfies the polynomial equation $\det(\BI-\BC_0\BT(t_1,t_2))=0$. The equation itself is a polynomial
equation in $\mathbb{R}$ and the boundary of $\CT$ is easily identifiable as a connected component of the solution set bounding a domain containing the origin and
no other solutions of the equation.
There then exists a vector $\Bv$ such that
\beq \Bv=\BC_0\BT\Bv. \eeq{2.3}
in which $\BT\Bv$ does not necessarily correspond to a symmetric matrix.
Here $\Bv$ represents a stress and in the basis \eq{1.1} has the representation
\beq \Bv=\bpm v_1 \\ v_2 \\ v_3 \\ v_4 \epm, \eeq{2.4}
where $v_3=v_4$ because the range of $\BC_0$ consists of symmetric matrices. Assuming $v_3=v_4$ is nonzero we may normalize $\Bv$ so that $v_3=v_4=1$.

Consider adding to $\BS_0-\BT$ a perturbing 
\beq \BGD=\bpm -\Ge_1 & 0 & 0 & 0 \\ 0 & -\Ge_2 & 0 & 0 \\ 0 & 0 & \Ge_1 & 0\\ 0 & 0 &0 & \Ge_2 \epm
\eeq{2.5}
to form a new matrix $\BS_0-\BT+\GD$. Note that the quadratic form associated with $\BT-\BGD$ is a null-Lagrangian and we are in effect
just perturbing $t_1$ and $t_2$. The corresponding change in the lowest eigenvalue is $\Bv^\dagger\BGD\Bv$ and the new bound becomes
$1/(2\Gm_*)\geq t_1+\Ge_1$ which will be better than the old bound if $\Ge_1>0$ and is valid for small $\Ge_1>0$ provided $\Bv^\dagger\BGD\Bv>0$, i.e. provided
$\BS_0-\BT+\BGD$ remains positive definite. To avoid a contradiction we must have that for all small $\Ge_1>0$ and $\Ge_2$ (not
necessarily positive),
\beq 0\leq -\Bv^\dagger\BGD\Bv=\Ge_1(|v_1|^2-|v_3|^2)+\Ge_2((|v_2|^2-|v_4|^2). \eeq{2.6}
Assuming $v_3=v_4=1$ this implies $|v_1|>1$ and $|v_2|=1$. We are free to rotate $\Bv$, by applying $\BR_\Gf$ to it, so that $e^{-i2\Gf}v_2=-1$. Accordingly, we replace
$\BS_0$ by the compliance tensor of the rotated crystal, which we redefine as our new $\BS_0$ that then has an associated value $v_2=-1$ (with $v_3=v_4=1$).

The bound $\BS_*-\BT\geq 0$ will be attained if we can find a divergence free stress field $\BGs(\Bx)$ that in the basis \eq{1.1} takes the form
\beq \BGs(\Bx)=\Ga(\Bx)\BR_{\Gt(\Bx)}\Bv, \eeq{2.7}
and is such that $\lang\BGs\rang\ne 0$
Then, with $\BC(\Bx)$ given by \eq{2.9}, the same argument as in \eq{2.34}-\eq{2.12} applies.

Second, to obtain a lower bound on $\Gm_*$ we require that $t_1$ and $t_2$ be chosen such that
\beq \BC_0-\BT\geq 0,\quad \det(\BC-\BT)=0. \eeq{2.13}
Among such pairs $(t_1,t_2)$ satisfying this we select the pair having the maximal value of $t_1$, since the bounds $\BC_*-\BT\geq 0$ imply
\beq 2\Gm_*\geq t_1,\quad 2\Gm_*\geq t_2,\quad \Gk_*(t_1+t_2)+t_1t_2\geq 0. \eeq{2.13a}
Defining this maximum value of $t_1$ to be $2\Gm_*^-$, we then have the lower bound $\Gm_*\geq\Gm_*^-$.
Finding $\Gm_*^-$ is accomplished by taking, in the $(t_1,t_2)$ plane, the simply connected loop of the curve $\det(\BC_0-\BT)=0$ that is closest to the origin, and surrounds it, and
finding the maximal value of $t_1$ along this curve. There then exists a vector $\Bw$ such that
\beq \BC_0\Bw=\BT\Bw, \eeq{2.14}
in which $\BT\Bw$ must be a symmetric matrix.
Here $\Bw$ represents a displacement gradient and in the basis \eq{1.1} has the representation
\beq \Bw=\bpm w_1 \\ w_2 \\ w_3 \\ w_4 \epm, \eeq{2.15}
where $t_1w_3=t_2w_4$ because  $\BT\Bw$ must be a symmetric matrix. Assuming $w_4$ is nonzero we may normalize $\Bw$ so that $w_3=t_2/t_1$ and $w_4=1$.

By adding to $\BC_0-\BT$ a perturbing matrix
\beq \BGD=\bpm -\Ge_1 & 0 & 0 & 0 \\ 0 & -\Ge_2 & 0 & 0 \\ 0 & 0 & \Ge_1 & 0\\ 0 & 0 &0 & \Ge_2 \epm
\eeq{2.16}
to form a new matrix $\BC_0-\BT+\GD$. Note that the quadratic form associated with $\BT-\BGD$ is a null-Lagrangian and we are in effect
just perturbing $t_1$ and $t_2$. The corresponding change in the lowest eigenvalue is $\Bw^\dagger\BGD\Bw$ and the new bound becomes
$\Gm_*>(t_1+\Ge_1)/2$ which will be better than the old bound if $\Ge_1>0$ and is valid for small $\Ge_1>0$ provided $\Bw^\dagger\BGD\Bw>0$, i.e. provided
$\BC_0-\BT+\BGD$ remains positive definite. To avoid a contradiction we must have that for all small $\Ge_1>0$ and $\Ge_2$ (not
necessarily positive),
\beq 0\leq -\Bw^\dagger\BGD\Bw=\Ge_1(|w_1|^2-|w_3|^2)+\Ge_2((|w_2|^2-|w_4|^2). \eeq{2.17}
Assuming $w_3=t_2/t_1$ and $w_4=1$ this implies $|w_1|>|t_2/t_1|$ and $|w_2|=1$. We are free to rotate $\Bw$, by applying $\BR_\Gf$ to it, so that $e^{-i2\Gf}w_2=1$. Accordingly, we replace
$\BS_0$ by the compliance tensor of the rotated crystal, which we redefine as our new $\BS_0$ that then has an associated value $w_2=1$ (with $w_3=t_2/t_1$ and $w_4=1$).

The bounds will be attained if we can find a displacement gradient $\BE(\Bx)$ that in the basis \eq{1.1} takes the form
\beq \BE(\Bx)=\Gb(\Bx)\BR_{\Gt(\Bx)}\Bw, \eeq{2.18}
and is such that the symmetric part of $\lang\BE\rang$ is nonzero. 
Then by the properties of $\BT$ and $\Bw$, $\BGs(\Bx)=\BT\BE(\Bx)$ is a divergence free and symmetric, i.e. it represents a stress field. Also we have
\beq \BGs(\Bx)=\BT\BE(\Bx)=\Gb(\Bx)\BR_{\Gt(\Bx)}\BT\Bw. \eeq{2.19}
Taking a polycrystal with in the basis \eq{1.1} the elasticity tensor \eq{2.20},
we see, using \eq{2.14}, that
\beq \BC(\Bx)\BE(\Bx)=\Gb(\Bx)\BC(\Bx)\BR_{\Gt(\Bx)}\Bw=\Gb(\Bx)\BR_{\Gt(\Bx)}\BC_0\Bw=\Gb(\Bx)\BR_{\Gt(\Bx)}\BT\Bw=\BGs(\Bx).
\eeq{2.21}
As the fields $\BE(\Bx)$ and $\BGs(\Bx)$ solve the constitutive relation, and the differential constraints, their averages are
related by the effective elasticity tensor $\BC_*$:
\beq \lang\BGs\rang=\BC_*\lang\BE\rang, \eeq{2.22}
and taking averages of both sides of $\BGs(\Bx)=\BT\BE(\Bx)$ gives
\beq \lang\BGs\rang=\BT\lang\BE\rang=\BC_*\lang\BE\rang. \eeq{2.23}
So the bound $\BC_*-\BT\geq 0$ is attained.

 %%%%%%%%%
 %%%%%%%%%%%%%%%%%%%%%%%%%%%%%%%%%%%%%%%%%%%%%%%%%%%%%%%%%%%%%%%%%%%%%%%%
\section{Proof of the simultaneous attainability of the lower bulk modulus and upper shear modulus bounds (point A)}
\setcounter{equation}{0}
%%%%%%%%%%%%%%%%%%%%%%%%%%%%%%%%%%%%%%%%%%%%%%%%%%%%%%%%%%%%%%%%%%%%%%%%%%
We look for geometries having a possibly anisotropic effective compliance tensor $\BS'$ such that for some $\Bv'\ne 0$ and $\Be'$,
each representing symmetric matrices,
\beq \BS'-\BT\geq 0,\quad \Bv'=\BC'\BT\Bv', \quad \BI:\Be'=k,\quad \BC'\Be'=\BI. \eeq{3.1}
Ultimately $\BC'$ will represent a positive definite effective tensor of a polycrystal with orthotropic symmetry
obtained from our starting material with elasticity tensor $\BC_0$ and following the construction in \cite{Avellaneda:1996:CCP}, it is then easy to go from $\BC'$ to a
polycrystal geometry with an isotropic effective tensor $\BC_*$ attaining the point A.
From these definitions we deduce that
\beq \overline{\Be'}:\Bv=\overline{\Be'}:\BC'\BT\Bv'=\overline{\BC'\Be'}:\BT\Bv'=\BI:\BT\Bv',
\eeq{3.2}
or equivalently, by normalization and by rotating and redefining $\BS'$ as necessary so that
\beq |v_1'|\geq 1, \quad v_2'=-1, \quad v_3'=v_4'=1\quad \text{and} \quad e_2'=\overline{e_1'},\quad e_3'=e_4'=k/2, \eeq{3.3}
\eq{3.2} reduces to
 \beq \overline{e'_1}v_1-e'_1=-\Gl,\quad\text{where}\quad \Gl=k+t_1+t_2. \eeq{3.4}
The last inequality in \eq{2.2a}, with $1/\Gk_*=k$ (so that the bulk modulus is attained), implies $\Gl\geq 0$. 

We will be studying trajectories that can be visualized as paths in the complex $e_1'(p)$-plane
parameterized by a real variable $p$.
Of particular interest are those trajectories that have loops with a tail from the $e_1$ associated with $\BC_0$ to the loop. As one goes around the trajectory
tail and loop there is an associated elasticity tensor $\BC'(e_1')$ having a unique value, modulo rotations, at the self intersection point of the trajectory. By introducing
a ``mirror'' material to $\BC_0$ we will see later in this section how to go from $e_1'$ to this associated tensor $\BC'(e_1')$ through hierarchical laminations
of the type illustrated in \fig{1}. 

In terms of $e_1'$ \eq{3.4} implies
\beq v_1'=(e_1'-\Gl)/\overline{e_1'}, \eeq{3.5}
and the constraint that $|v_1'|\geq 1$ holds if and only if
\beq \Real(e_1')\leq \Gl/2.
\eeq{3.6}

Now let us consider the case that $e_1'$ is real. Then  $v_1'=(e_1'-\Gl)/e_1'$ is also real.
The real parts of the matrices corresponding to $\Be'$, $\BC'\Be'=\BI$, $\Bv'=\BC'\BT\Bv'$ and $\BT\Bv'$ are all
diagonal matrices. This implies $\BC'$ is orthotropic since it maps any diagonal real matrix, being a linear combination of the real parts of the matrices
corresponding to $\Be'$ and $\BT\Bv'$, to a diagonal real matrix. We consider a layering in direction $\Bn=(0,1)$ of the two stress matrices
\beq \BGs_1=e^{i\Gt}\BR_{\Gt}\Bv',\quad \text{and}\quad \BGs_2=e^{-i\Gt}\BR_{-\Gt}\Bv'. \eeq{3.10}
These will be compatible if their difference
\beq \BGs_1-\BGs_2=\bpm 2iv_1'\sin{3\Gt} \\ 2i\sin{\Gt} \\ 2i\sin{\Gt}\\ 2i\sin{\Gt} \epm \eeq{3.11}
is of the form \eq{1.12}. Assuming $\Gt\ne 0$, this will be the case if
\beq v_1'=\frac{\sin\Gt}{\sin(3\Gt)}=\frac{1}{2y+1},\quad\text{where}\quad y=\cos(2\Gt), \eeq{3.12}
giving
\beq e_1'=\frac{\Gl}{1-v_1'}=\frac{\Gl(2y+1)}{2y}.\eeq{3.12a}
Note that, as $\Gt$ is real, $y=\cos{2\Gt}$ must lie between $1$ and $-1$. In fact we need
\beq -1\leq y \leq 0, \eeq{3.12b}
to ensure that the constraint $|v_1|\geq 1$ is satisfied.
The average field will then be
\beq \Bv''=p\BGs_1+(1-p)\BGs_2=\bpm v_1'\cos{3\Gt} \\ \cos{\Gt} \\ \cos{\Gt}\\ \cos{\Gt} \epm +(2p-1)\bpm i\sin{\Gt} \\ i\sin{\Gt} \\ i\sin{\Gt}\\ i\sin{\Gt} \epm.
\eeq{3.13}
At the same time, in connection with the bulk modulus Reuss-Hill bound, we layer in direction $\Bn=(0,1)$ the two strain matrices
\beq \BGe_1=\BR_{\Gt}\Be',\quad \text{and}\quad \BGe_2=\BR_{-\Gt}\Be',
\eeq{3.14}
which are compatible since $\Bt\cdot(\BGe_1-\BGe_2)\Bt=0$, with $\Bt=(1,0)$. Note that as we are focusing on the Reuss-Hill bound the corresponding
stress matrices $\BGs_1$ and $\BGs_2$ will both equal the identity matrix, as implied by \eq{2.30b}. Such matrices when expressed in the basis
\eq{1.1} cannot be normalized as their second element is zero.
This produces the average displacement gradient
\beq \Be''=p\BGe_1+(1-p)\BGe_2=\bpm e_1'\cos(2\Gt) \\ e_1'\cos(2\Gt) \\ k/2 \\ k/2 \epm+(p-0.5)\bpm ie_1'\sin(2\Gt) \\ -ie_1'\sin(2\Gt) \\ 0 \\ 0 \epm.
  \eeq{3.15}
  We next normalize and rotate the average fields, using the rotation
  \beq \BR=\bpm  (1-i\Gj)/(1+i\Gj) & 0 & 0 & 0 \\ 0 & (1+i\Gj)/(1-i\Gj) & 0 & 0 \\ 0 & 0 & 1 & 0\\ 0 & 0 &0 & 1 \epm, \quad\text{with}\quad\Gj=(2p-1)\tan\Gt,
  \eeq{3.16}
to obtain
\beq \Bv=\frac{\BR\Bv''}{(1+i\Gj)\cos{\Gt}}=\bpm v_1 \\ -1 \\ 1 \\ 1 \epm,\quad \Be=\BR\Be''=\bpm e_1 \\ \overline{e_1} \\ k/2 \\ k/2 \epm,
\eeq{3.17}
with
\beq v_1=\frac{\left(\frac{\tan\Gt}{\tan(3\Gt)}+i\Gj\right)(1-i\Gj)}{(1+i\Gj)^2},\quad
e_1=-\frac{e_1'\left(\cos(2\Gt)+i\Gj\frac{\sin(2\Gt)}{\tan\Gt}\right)(1-i\Gj)}{(1+i\Gj)}. \quad
\eeq{3.18}
One can double check that the identity $\overline{e_1}v_1-e_1=-\Gl$ holds, as it should. In terms of $y=\cos{2\Gt}$ and $\Gj$ the expression for $e_1$ becomes
\beq e_1=e_1(y,\Gj)\equiv\frac{\Gl(2y+1)\left(1+i\Gj(1+\frac{1}{y})\right)(1-i\Gj)}{2(1+i\Gj)}, \eeq{3.19}
where we have used \eq{3.12a}.

Values of $\Gj$ such that $\Gj^2\leq (1-y)/(1+y)=(\tan\Gt)^2$ which thus correspond to values of $p$ between $0$ and $1$ give an optimal composite that is polycrystal laminate of an orthotropic material. It is built from two reflected orientations of the orthotropic material, reflected about the layering direction. This is not so interesting as our objective is to build
an optimal polycrystal built from a non-orthotropic crystal. On the other hand, as we will see, values of $\Gj$ such that $\Gj^2\geq (1-y)/(1+y)$ mean that we
can optimally layer a material that is the original crystal with an associated
$e_1$ with a rotated orthotropic material to obtain a composite having an effective tensor that is a different rotation
of the same orthotropic material. The volume fraction occupied
by the original crystal in this laminate is $q=1/(1-p)$ if $p<0$ or $q=1/p$ if $p>1$. We now just treat the case where $p<0$ as the case $p>1$ can be treated
similarly. 

In this laminate we can repeatedly replace the rotated orthotropic material with a rotated laminate of the original crystal
and the same orthotropic material in an appropriate orientation, as illustrated in \fig{1}, where the blue material represents the orthotropic material
and now \fig{1}(a) has the same volume fractions as 1(b), 1(c), and so forth.
\fig{1}(a) has itself the same effective tensor as the blue material and the geometry is fully self-similar.
The values of $e_1$, $v_1$, $e_1'$, $v_1'$, $p$ and $\Gt$ 
give us stress and strain fields that solve the elasticity equations at each stage.
Doing this ad infinitum so that the orthotropic material occupies a vanishingly small volume fraction,
we obtain a polycrystal of the material  corresponding to the value $e_1$ that has an orthotropic effective tensor corresponding to the value $e_1'$. Then at the final stage we
can construct an optimal elastically isotropic polycrystal from the orthotropic polycrystal corresponding to $e_1'$ using, for example, the prescription outlined in
\cite{Avellaneda:1996:CCP}.

The question is now: what is the range of values taken by $e_1$ as $y$ and $\Gj$ range over the set
\beq \BGS=\{(y,\Gj)~|~-1<y<0,\quad \Gj^2\geq (1-y)/(1+y)\}? \eeq{3.20}
The line $y=-1/2$, parameterized by $\Gj$, is rather singular as the whole line gets mapped to $e_1=0$. So rather than considering the image of the
set $\BGS$, let us consider the image of the two sets
\beqa  \BGS_1 & = & \{(y,\Gj)~|~-1<y<-1/2,\quad \Gj^2\geq (1-y)/(1+y)\},\nonum
\BGS_2 & = & \{(y,\Gj)~|~-1/2<y<0,\quad \Gj^2\geq (1-y)/(1+y)\}.
\eeqa{3.21}
Our aim is to show that $e_1(y,\Gj)$ maps $\BGS_1\cup\BGS_2$ to the set
\beq \BGO=\{e_1~|~\Real(e_1)\leq \Gl/2, \Real(e_1)\ne 0, \Imag(e_1)\ne 0\}, \eeq{3.21a}
in a one-to-one and onto (bijective) fashion. 

From \eq{3.19} we have
\beqa \Real(e_1) & = & \frac{\Gl(2y+1)\left(1+\Gj^2+\frac{2\Gj^2}{y}\right)}{2(1+\Gj^2)}, \nonum
 \Imag(e_1) & = & -\frac{\Gl(2y+1)\Gj[\Gj^2(1+y)-(1-y)]}{2y(1+\Gj^2)} .
\eeqa{3.22}
First we look at what happens when $(y,\Gj)$ is close to the boundary of the closure of $\BGS_1$ or the boundary of the closure of $\BGS_2$. To start, consider those points
in $\BGS_1$ and $\BGS_2$ near the line $y=-1/2$ when $|\Gj|$ is very large. From \eq{3.22} it follows that
\beq \Real(e_1) \approx -3\Gl(2y+1)/2,\quad \Imag(e_1)\approx \Gl q \Gj, \quad\text{where}\quad q=(2y+1)\Gj/2. \eeq{3.23}
Thus $\Imag(e_1)$ can take any nonzero value with $\Real(e_1)$ arbitrarily close to $0$. In both regions, near the curves $\Gj^2\geq (1-y)/(1+y)$ one has that $\Imag(e_1)$ is close to
zero and
\beq \Real(e_1) \approx \Gl(2y+1)/(2y), \eeq{3.24}
which varies from $\Gl/2$ to minus infinity as $y$ varies from $-1$ to $0$. Again in both regions, but when $|\Gj|$ is large and $y$ is not close to $-1/2$,
\beq \Real(e_1)\approx \Gl(2y+1)\left(1+\frac{2}{y}\right)/2,\quad \Imag(e_1) \approx -\frac{\Gl(2y+1)\Gj[(1+y)]}{2y}. \eeq{3.25}
Thus $\Imag(e_1)$ is arbitrarily large and $\Real(e_1)$ ranges from $\Gl/2$ to minus infinity as $y$ varies from $-1$ to $0$. 
In region $\BGS_1$ when $\Gj$ is large and $y$ is close to $-1$, suppose that $\Gj^2-(1-y)/(1+y)$ scales in proportion to $\Gj$, i.e. $\Gj^2-(1-y)/(1+y)=h|\Gj|$, with
for some $h>0$. Then from \eq{3.22} it follows that
\beq \Real(e_1) \approx \Gl/2,\quad \Imag(e_1)\approx -\Gl h\Gj/(2|\Gj|).
\eeq{3.26}
So $\Real(e_1)$ can be arbitrarily close to the line $\Real(e_1)\approx \Gl/2$ with an arbitrary nonzero value of $ \Imag(e_1)$. 
Finally, we look at what happens in $\BGS_2$ near the curve $\Gj^2\geq (1-y)/(1+y)$ when $y$ is small (and hence $\Gj^2$ is close to $1$). Now \eq{3.22} implies
\beq \Real(e_1)\approx \frac{\Gl}{y}, \quad \Imag(e_1)\approx -\Gl\Gj r \quad\text{where}\quad r=-[\Gj^2(1+y)-(1-y)]/(2y). \eeq{3.26a}
In other words, $\Real(e_1)$ can be arbitrarily large, while $\Imag(e_1)$ can take any desired value.

To check if the mapping is locally one-to-one (injective), we look at the derivatives,
\beqa \frac{de_1}{dy} & = & \frac{\Gl[2-i\Gj(1/y^2-2)](1-i\Gj)}{2(1+i\Gj)}, \nonum
\frac{de_1}{d\Gj} & = & \frac{\Gl(2y+1)\{2\Gj(1+1/y)-i[(1-1/y)-\Gj^2(1+1/y)]\}}{2(1+i\Gj)^2},
\eeqa{3.27}
and their arguments,
\beqa \arg\left(\frac{de_1}{dy}\right)& = & \arg\left[2-i\Gj(1/y^2-2)\right]+2\arg(1-i\Gj), \nonum
\arg\left(\frac{de_1}{d\Gj}\right)& = & \arg\left\{2-i\Gj\left[\frac{y-1}{\Gj^2(y+1)}-1\right]\right\}+2\arg(1-i\Gj).
\eeqa{3.28}
Clearly both derivatives in \eq{3.27} are nonzero in the regions $\BGS_1$ and $\BGS_2$ (but $de_1/d\Gj$ is zero along the
line $y=-1/2$). Also, from the inequality
\beq (1/y^2-2)-\left[\frac{y-1}{\Gj^2(y+1)}-1\right]=(1/y^2-1)+\left[\frac{1-y}{\Gj^2(y+1)}\right]>0 \eeq{3.29}
we see that $\arg(de_1/dy)\ne\arg(de_1/d\Gj)$. Thus the mapping is locally one-to-one. This completes the proof that $e_1(y,\Gj)$ maps $\BGS_1\cup\BGS_2$ to $\GO$
in a one-to-one and onto fashion (i.e. it is a bijective mapping). Note also that $\GO$ includes all $e_1$ that correspond to positive definite $\BC_0$, excluding those  positive definite $\BC_0$ that have $\Real(e_1)=0$ or  $\Imag(e_1)=0$. 

In practice, given $e_1$ we can find $\Gj$ and $\Gt$ by solving
\beq |e_1|^2=\Gl^2(2y+1)^2\left[1+\Gj^2\left(1+\frac{1}{y}\right)^2\right],\quad \arg(e_1)=\tan^{-1}\left[\Gj\left(1+\frac{1}{y}\right)\right]-2\tan^{-1}(\Gj).
\eeq{3.50}
The first equation can be easily solved for $\Gj$ in terms of $y$, and then substituted in the second equation which can then be numerically solved
for $y$. One picks the value of $y$ satisfying \eq{3.12b} that has an associated value of $\Gj$ satisfying $\Gj^2\geq (1-y)/(1+y)$.

Now given a fourth-order elasticity tensor $\BC_0$ that is positive definite on the space of symmetric matrices, and which has the associated 
complex number $e_1'=e_1$, we need to show that given a trajectory having a tail from $e_1$ to a loop, there is an associated $\BC'(e'_1)$
on the tail and loop that is also positive definite on the space of symmetric matrices. To do this we introduce a ``mirror material''
with elasticity tensor $\BC_0^M$ obtained by reflecting $\BC_0$  about the $x_1$ axis. Under the reflection $x_1\to -x_1$ the basis
vectors $\BGr_1$, $\BGr_2$, $\BGr_3$, and $\BGr_4$ transform to $\BGr_2$, $\BGr_1$, $\BGr_4$, and $\BGr_3$ respectively, and $\BT(t_1,t_2)$ transforms to $\BT(t_2,t_1)$.

Hence the identity $\tilde{\Bv}=\BC_0^M\BT(t_2,t_1)\tilde{\Bv}$ will be satisfied by a matrix $\tilde{\Bv}$
that has the vector representation
\beq \tilde{\Bv}=\bpm v_2 \\ v_1 \\ v_4 \\ v_3 \epm, \eeq{3.51}
where the $v_i$ are the elements of the vector $\Bv$ satisfying $\Bv=\BC_0\BT\Bv$. Further, as $\BGr_2$, $\BGr_1$, $\BGr_4$, $\BGr_3$,
and $\BT(t_2,t_1)$ transform back to  $\BGr_1$, $\BGr_2$, $\BGr_3$, $\BGr_4$, and  $\BT(t_1,t_2)$ under complex conjugation, by
taking complex conjugates of the identity
\beq v_2\BGr_1+v_1\BGr_2+v_4\BGr_3+v_3\BGr_4=\BC_0^M\BT(t_2,t_1)(v_2\BGr_1+v_1\BGr_2+v_4\BGr_3+v_3\BGr_4), \eeq{3.52}
we see that $\Bv^M=\BC_0^M\BT\Bv^M$ is satisfied by the matrix $\Bv^M$ that has the representation
\beq \Bv^M=\bpm \overline{v}_1 \\ \overline{v}_2 \\ \overline{v}_3 \\ \overline{v}_4 \epm. \eeq{3.53}
Similarly, the vector $\Be$ gets replaced by its complex conjugate and hence $\BC_0^M$ is represented by $\overline{e}_1$. For an appropriate rotations
of $\BC_0$ in the laminate, with $\BC_M$ in the laminate being rotated by the same angle but in the reverse direction to maintain mirroring of the phases. 
the fields we have obtained solve the cell problem, and the three materials X,Y and Z in \fig{3} are represented by real values of $e'_1$ and
therefore orthotropic. Moreover, since an orthotropic tensor is invariant under mirroring about the axis of orthotropy,
the effective elasticity tensors X and Y, as they are clearly  mirror images of each other,
must also be rotations of each  other. We can express them as rotations by angles $\Gt$ and $-\Gt$ of an orthotropic tensor $\BC'$ with axes the same as the coordinate axes,
thus defining their effective elasticity tensors, $\BC_{\Gt}$ for \fig{3}(a)  and $\BC_{-\Gt}$ for \fig{3}(b).

We now just treat the case where $p<0$ as the case $p>1$ can be treated
similarly. Referring to \fig{4} we see that the volume fractions of $\BC_0$ in $\BC_{\Gt}$ and $\BC_{-\Gt}$ (in these laminates of $\BC_0$ and $\BC_0^M$)
are respectively $(1-p)/(1-2p)$ and $-p/(1-2p)$.
Note further from \fig{3} and \fig{4} that by laminating $\BC_0$ with $\BC_{-\Gt}$, again in the direction $\Bn=(0,1)$, with $\BC_0$ occupying a volume fraction $1/(1-p)$
one obtains $\BC_\Gt$. 
%Note that $\Be'$ is also real, and the mixture consequently is orthotropic, whe. However, the $\BC'$ we need corresponds to the self intersection point
%of the lamination trajectory and this requires that the prescribed unequal volume fractions of $\BC_0$ mixed with (b) to obtain (a) are $(1-p)/(1-2p)$ and $-p/(1-2p)$
%respectively.
\begin{figure}[ht!]
\includegraphics[width=0.9\textwidth]{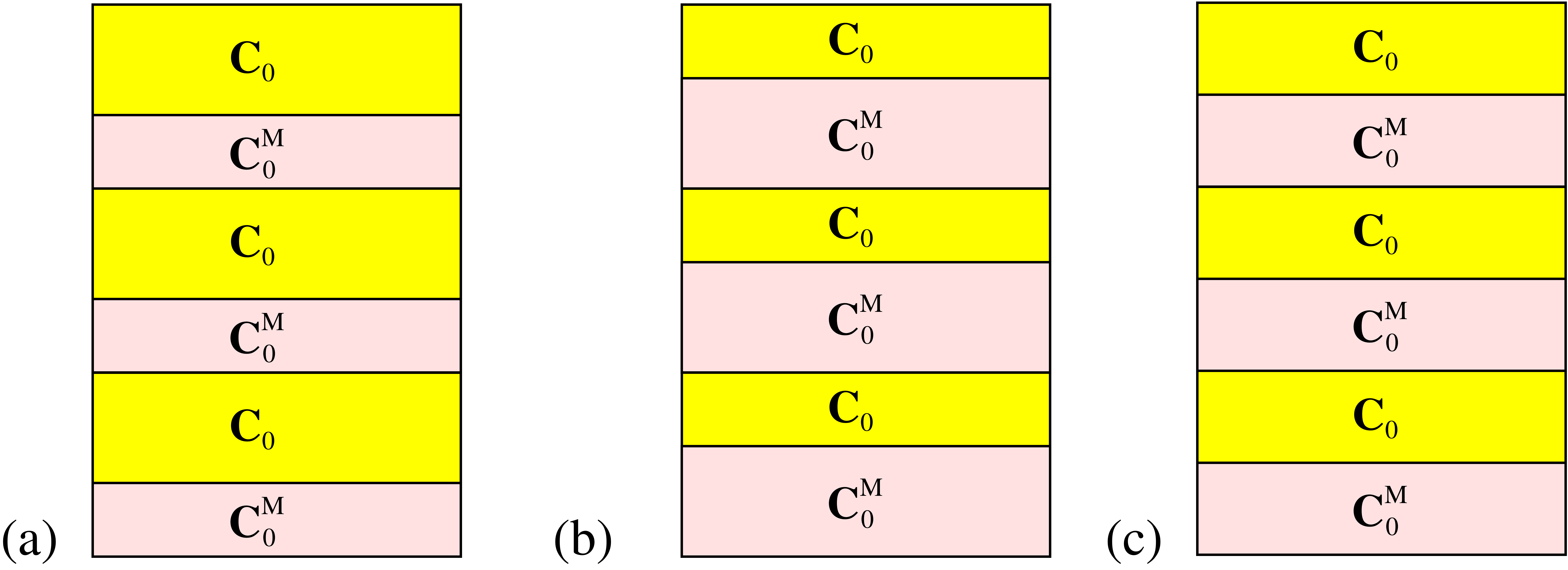}
\caption{Laminates of properly orientated tensors $\BC_0$ and $\BC_M$ that have fields with real values of $e_1'$, and hence which have orthotropic symmetry.
This orthotropy is immediately evident in Z less so in X and Y. 
Observe that by laminating $\BC_0$ with Y, again in the direction $\Bn=(0,1)$, one obtains X.
}
\labfig{3}
\end{figure}

\begin{figure}[ht!]
\includegraphics[width=0.9\textwidth]{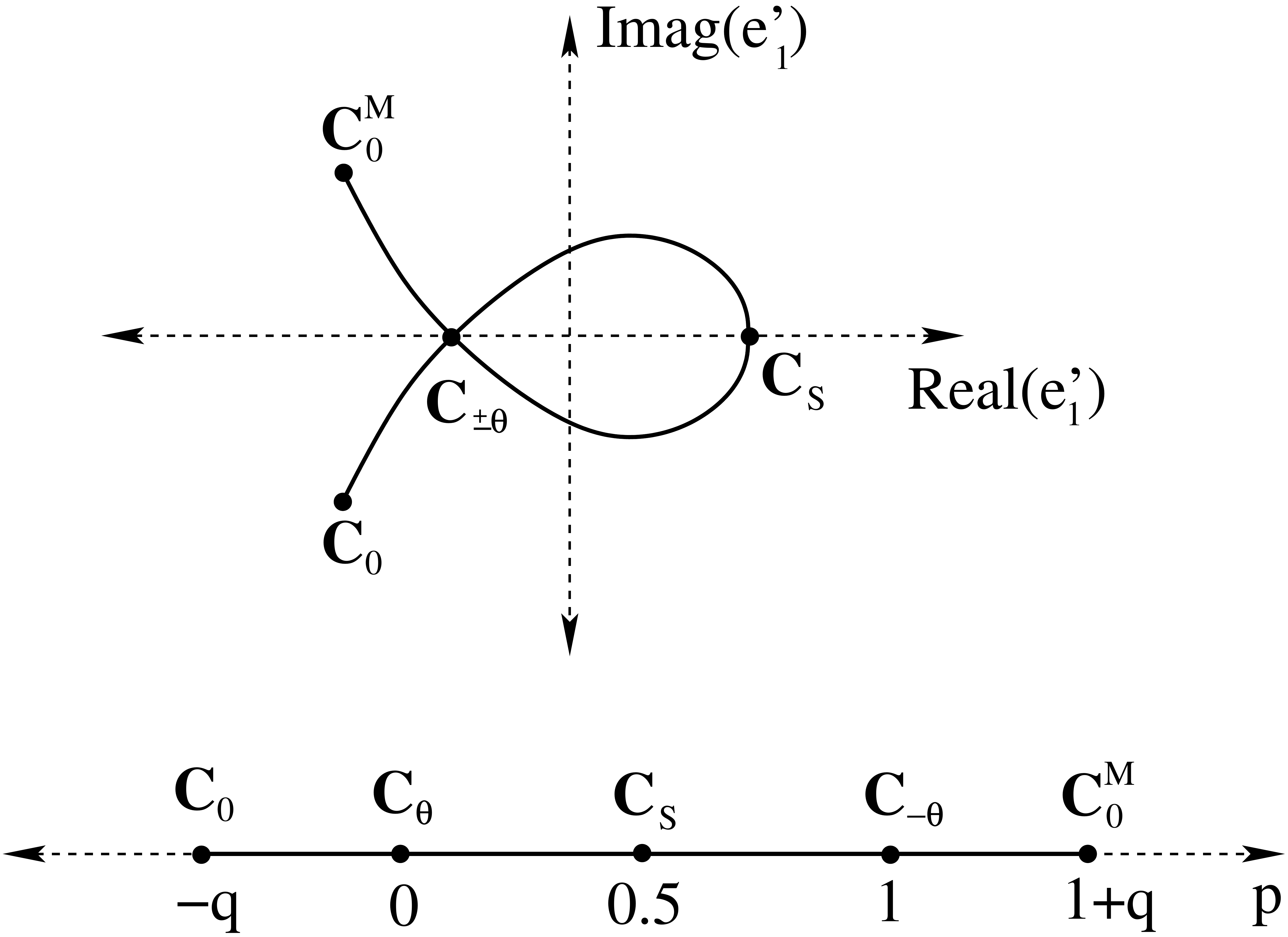}
\caption{Figure (a) shows a schematic representation of the trajectory in the $\Be_1$ plane. The self intersection point corresponds to the geometries X and Y
in Figure 3.  When parameterized by $p$ and straightened the trajectory becomes (b). As seen from (b) to obtain $\BC_{\Gt}$ one should mix $\BC_0$ with $\BC_0^M$ in proportions $(1-p)/(1-2p)$ and $-p/(1-2p)$ respectively, or alternatively mix  $\BC_0$ with $\BC_{-\Gt}$ in proportions $1/(1-p)$ and $-p/(1-p)$. In either case the proportion of $\BC_0$ in the mixture is greater than the  other material
}
\labfig{4}
\end{figure}

An explicit formula for $\BC'$ can be obtained using the lamination formula of Francfort and Murat \cite{Francfort:1986:HOB}, or more simply
from the formulae (5.22) and (5.26) in \cite{Milton:1992:CMP} which are easily generalized to allow for unequal volume fractions. 
As both  $\BC_0$ and $\BC_0^M$ are positive definite tensors on the space of symmetric matrices, so must be the orthotropic material $\BC'$. It may be the case
that the elasticity tensor $\BC_0^M$ is available as a constituent material. That happens when the rotations $\BR(\Bx)$ in \eq{1.b} are not restricted to be proper rotations, but can include
reflections as well (having $\det[\BR(\Bx)]=-1$).  Then we do not require an infinite rank lamination scheme to obtain the orthotropic material $\BC'$ represented by $\Be'$ as it is a laminate of $\BC_0$ and $\BC_0^M$. Given other values of $\Be'$ on the trajectory we can find the associated $\BC'$, not necessarily orthotropic,
from the lamination formula with the appropriate volume fraction in the laminate of $\BC_0$ and $\BC_0^M$.

The above argument avoids the difficult question as to what values of $e_1$ correspond to some positive semidefinite elasticity tensor $\BC$. Certainly
the positivity of $\Bv^\dagger\cdot\BS_0\Bv=\Bv^\dagger\cdot\BT\Bv$ with $|v_2|=v_3=v_4=1$ implies $|v_1|^2\geq 1$, and a tighter bound is
obtained from the inequality
\beq (\Bv+\Ga\BI)^\dagger\BS_0(\Bv+\Ga\BI)=t_1(|v_1|^2-1)+\Ga^2k+2\Real[\Ga(\overline{v_1}e_1+\overline{v_2}\overline{e_1}+k)],
\eeq{3.54}
which holds for all complex $\Ga$. It may be the case that even tighter bounds on $e_1$ are needed to guarantee that $e_1$ corresponds to some
positive semidefinite elasticity tensor $\BC$. Fortunately, we do not need them.

Having obtained an orthotropic material attaining the bounds that corresponds to the point $e_1'$ (with $e_1'$ being real) we can follow the procedure in
\cite{Avellaneda:1996:CCP} to obtain an isotropic material attaining the bounds. More generally, it is easily seen that any value of $e_1$ inside the loop that self intersects
at $e'_1$ is attainable.

It remains to treat the special cases when $e_1$ is purely real or $e_1$ is purely imaginary. The first corresponds to an orthotropic material
and is treated in \cite{Avellaneda:1996:CCP}. In the second case, using the fact that $\overline{e_1}=-e_1$ and $v_1=(e_1-\Gl)/\overline{e_1}=-1+\Gl/e_1$,
consider the fields
\beqa \BGs_1 & = & \BR_{\pi/2}\Bv=\bpm 1-\Gl/e_1 \\ 1 \\ 1 \\ 1 \epm,\quad \BGs_2=\bpm -\Gl/e_1 \\ 0 \\ 0 \\ 0 \epm, \nonum
\BGe_1 & = & \BR_{\pi/2}\Be=\bpm -e_1 \\ e_1 \\ k/2 \\ k/2 \epm,\quad \BGe_2=\bpm 0 \\ 0 \\ k/2 \\ k/2 \epm.
\eeqa{3.55}
As the differences
\beq \BGs_1-\BGs_2=\bpm 1 \\ 1 \\ 1 \\ 1 \epm,\quad \BGe_1 -\BGe_2 =\bpm -e_1 \\ e_1 \\ 0 \\ 0 \epm
\eeq{3.56}
are of the same form as in \eq{1.12}, the fields are compatible. Also $\BGs_2$ and $\BGe_2$ represent shear and hydrostatic fields in an
isotropic medium with effective shear and bulk moduli $\Gm_*^+$ and $\Gk_*^-$. So there exists a suitable trajectory joining
this isotropic effective medium with the original crystal. Using an infinite rank lamination scheme that corresponds to the effective
medium approximation, similar to the one detailed at the beginning of section 4.2 of \cite{Avellaneda:1996:CCP}, we conclude that the point A is attained
even in these special cases. Note that this infinite rank lamination scheme is a limiting case of the infinite rank lamination scheme that attains points in $\GO$, in which the loop in the
trajectory shrinks to zero. 

 %%%%%%%%%%%%%%%%%%%%%%%%%%%%%%%%%%%%%%%%%%%%%%%%%%%%%%%%%%%%%%%%%%%%%%%%
\section{Proof of the simultaneous attainability of the lower bulk modulus and lower shear modulus bounds (point D)}
\setcounter{equation}{0}
%%%%%%%%%%%%%%%%%%%%%%%%%%%%%%%%%%%%%%%%%%%%%%%%%%%%%%%%%%%%%%%%%%%%%%%%%%
We look for geometries having a possibly anisotropic effective elasticity tensor $\BC'$ such that for some $\Bw'\ne 0$ and $\Be'$,
the latter representing a symmetric matrix,
\beq \BC'-\BT\geq 0,\quad \BC'\Bw'=\BT\Bw', \quad \BI:\Be'=k,\quad \BC'\Be'=\BI. \eeq{4.1}
Ultimately $\BC'$ will represent an effective
tensor of a polycrystal obtained from our starting material with elasticity tensor $\BC_0$. As $\BT\Bw'=\BC'\Bw'$ is a symmetric matrix, we have
that
\beq t_1w_3'=t_2w_4'. \eeq{4.2}
Also the definitions \eq{4.1} imply
\beq \overline{\Be'}:\BT\Bw=\overline{\Be'}:\BC'\Bw'=\overline{\BC'\Be'}:\Bw'=\BI:\Bw',
\eeq{4.3}
or equivalently, by normalization and by rotating and redefining $\BC'$ as necessary so that
\beq |w_1'|\geq |t_2/t_1|, \quad w_2'=1, \quad w_3'=t_2/t_1, \quad w_4=1 \quad \text{and} \quad e_2'=\overline{e_1'},\quad e_3'=e_4'=k/2, \eeq{4.4}
\eq{4.3} reduces to
\beq t_1\overline{e'_1}v_1+t_2e'_1=t_2\Gn,\quad\text{where}\quad \Gn=k+1/t_1+1/t_2. \eeq{4.5}
Note that the last inequality in \eq{2.13a}, with $k=1/\Gk_*$ (so that the bulk modulus is attained), implies $t_1t_2\Gn\geq 0$.

We can think of any material attaining the bounds as being parameterized by the complex number $e_1'$. In terms of it \eq{3.4} implies
\beq w_1'=(\Gn-e_1')t_2/(t_1\overline{e_1'}), \eeq{4.6}
and the constraint that $|w_1'|\geq |w_3'|=|t_2/t_1|$ holds if and only if
\beq \Gn(\Real(e_1')-\Gn/2)\geq 0.
\eeq{4.7}
Now let us further assume that $\BC'$ has orthotropic symmetry which implies that $e_1'$ is real
(and hence $w_1'=(\Gn-e_1')t_2/(t_1e_1')$ is also real).
We consider a layering in direction $\Bn=(0,1)$ of the two displacement gradients
\beq \BE_1=e^{i\Gt}\BR_{\Gt}\Bw',\quad \text{and}\quad \BE_1=e^{-i\Gt}\BR_{-\Gt}\Bw'. \eeq{4.8}
These will be compatible if their difference
\beq \BE_1-\BE_2=\bpm 2iw_1'\sin{3\Gt} \\ -2i\sin{\Gt} \\ 2it_2\sin{\Gt}/t_1 \\ 2i\sin{\Gt} \epm \eeq{4.9}
is of the form \eq{1.12}. Assuming $\Gt\ne 0$, this will be the case if
\beq w_1'=-\frac{t_2\sin\Gt}{t_1\sin(3\Gt)}=-\frac{t_2}{t_1(2y+1)},\quad\text{where}\quad y=\cos(2\Gt), \eeq{4.10}
giving
\beq e_1'=\frac{\Gn}{1+t_1w_1'/t_2}=\frac{\Gn(2y+1)}{2y}.\eeq{4.11}
Note that as $\Gt$ is real $y=\cos{2\Gt}$ must lie between $1$ and $-1$ and in fact we need
\beq -1\leq y \leq 0, \eeq{4.12}
to ensure that the constraint $|w_1|\geq |t_2/t_1|$ is satisfied.
The average field will then be
\beq \Bw''=p\BE_1+(1-p)\BE_2=\bpm w_1'\cos{3\Gt} \\ \cos{\Gt} \\ t_2\cos{\Gt}/t_1\\ \cos{\Gt} \epm +(2p-1)\bpm -it_2\sin{\Gt}/t_1 \\ -i\sin{\Gt} \\ it_2\sin{\Gt}/t_1\\ i\sin{\Gt} \epm.
\eeq{4.13}
At the same time, in connection with the bulk modulus Voigt bounds, we layer in direction $\Bn=(0,1)$ the two compatible strain matrices
\beq \BGe_1=\BR_{\Gt}\Be',\quad \text{and}\quad \BGe_2=\BR_{-\Gt}\Be',
\eeq{4.14}
to produce the average strain
\beq \Be''=p\BGe_1+(1-p)\BGe_2=\bpm e_1'\cos(2\Gt) \\ e_1'\cos(2\Gt) \\ k/2 \\ k/2 \epm+(p-0.5)\bpm ie_1'\sin(2\Gt) \\ -ie_1'\sin(2\Gt) \\ 0 \\ 0 \epm.
  \eeq{4.15}
  We next normalize and rotate the average fields, using the rotation
  \beq \BR=\bpm  (1-i\Gj)/(1+i\Gj) & 0 & 0 & 0 \\ 0 & (1+i\Gj)/(1-i\Gj) & 0 & 0 \\ 0 & 0 & 1 & 0\\ 0 & 0 &0 & 1 \epm, \quad\text{with}\quad\Gj=(2p-1)\tan\Gt,
  \eeq{4.16}
to obtain
\beq \Bw=\frac{\BR\Bw''}{(1+i\Gj)\cos{\Gt}}=\bpm w_1 \\ 1 \\ t_2/t_1 \\ 1 \epm,\quad \Be=\BR\Be''=\bpm e_1 \\ \overline{e_1} \\ k/2 \\ k/2 \epm,
\eeq{4.17}
with
\beq w_1=-\frac{t_2\left(\frac{\tan\Gt}{\tan(3\Gt)}+i\Gj\right)(1-i\Gj)}{t_1(1+i\Gj)^2},\quad
e_1=\frac{e_1'\left(\cos(2\Gt)+i\Gj\frac{\sin(2\Gt)}{\tan\Gt}\right)(1-i\Gj)}{(1+i\Gj)}. \quad
\eeq{4.18}
One can double check that the identity $t_1\overline{e_1}w_1/t_2+e_1=\Gn$ holds, as it should. In terms of $y=\cos{2\Gt}$ and $\Gj$ the expression for $e_1$ becomes
\beq e_1=e_1(y,\Gj)\equiv\frac{\Gn(2y+1)\left(1+i\Gj(1+\frac{1}{y})\right)(1-i\Gj)}{2(1+i\Gj)}, \eeq{4.19}
where we have used \eq{4.11}. Comparing \eq{3.6} with \eq{4.7} and \eq{4.19} with \eq{3.19} we see we can exactly repeat the subsequent analysis in the previous section
(with $\Gl$ replaced by $-\Gn$) to establish that the
lower bulk modulus and lower shear modulus bounds can be simultaneously attained.
 %%%%%%%%%%%%%%%%%%%%%%%%%%%%%%%%%%%%%%%%%%%%%%%%%%%%%%%%%%%%%%%%%%%%%%%%
\section{Microstructures that simultaneously attain the upper bulk modulus and upper shear modulus bounds (point B)}
\setcounter{equation}{0}
%%%%%%%%%%%%%%%%%%%%%%%%%%%%%%%%%%%%%%%%%%%%%%%%%%%%%%%%%%%%%%%%%%%%%%%%%%
We look for geometries having a possibly anisotropic effective compliance tensor $\BS'$ such that for some $\Bv'\ne 0$ and $\Bc'\ne 0$,
each representing symmetric matrices,
\beq \BS'-\BT\geq 0,\quad \BS'-\BT_0\geq 0,\quad\Bv'=\BC'\BT\Bv', \quad\Bc'=\BC'\BT_0\Bc'. \eeq{5.1}
Ultimately $\BC'$ will represent an effective
tensor of a polycrystal obtained from our starting material with elasticity tensor $\BC_0$, but in this section and the next one $\BC'$
will not typically have orthotropic symmetry. From these definitions we
deduce that
\beq \overline{\Bc'}:\BT\Bv=\overline{\BC'\BT_0\Bc'}:\BT\Bv'=\overline{\Bc'}:\BT_0\Bv',
\eeq{5.2}
or equivalently, by normalization and by rotating and redefining $\BS'$ as necessary so that
\beq |v_1'|\geq 1, \quad v_2'=1, \quad v_3'=v_4'=1\quad \text{and} \quad c_2'=\overline{c_1'},\quad c_3'=c_4'=1, \eeq{5.3}
\eq{5.2} reduces to
\beq \Ga_1\overline{c'_1}v'_1+\Ga_2c'_1=1, \eeq{5.4}
where
\beq \Ga_1=(t_0+t_1)/(2t_0+t_1+t_2),\quad\Ga_2=(t_0+t_2)/(2t_0+t_1+t_2)=1-\Ga_1 . \eeq{5.4a}
Note that the last inequality in \eq{2.2a}, with $1/\Gk_*=2t_0$ (so that the bulk modulus is attained), implies that $2t_0+t_1+t_2\geq 0$.
Also, since $t_0$ and $t_1$ are both positive, we conclude that $\Ga_1\geq 0$.

We can think of any material attaining the bounds as being parameterized by the complex number $c_1'$. In terms of it \eq{5.4} implies
\beq v_1'=\frac{1-\Ga_2c'_1}{\Ga_1\overline{c_1'}}, \eeq{5.5}
and the constraint that $|v_1'|\geq 1$ holds if and only if
\beq |1-\Ga_2c_1'|^2\geq\Ga_1^2|c_1'|^2.
\eeq{5.6}
This is automatically satisfied if $|c_1'|\leq 1$.

Now consider the stress field trajectory
\beq \BGs(\Gj)=\bpm e^{i2\Gt}v_1+e^{-i\Gt}\Gj \\  e^{-i2\Gt}+e^{-i\Gt}\Gj  \\ 1+e^{-i\Gt}\Gj \\ 1+e^{-i\Gt}\Gj \epm,
\eeq{5.7}
that we will associate with the upper shear modulus bounds, and the stress field trajectory
\beq \BGs_0(\Gj)=\bpm e^{i2\Gt}c_1+\Gf\Gj \\  e^{-i2\Gt}\overline{c_1}+\Gf\Gj  \\ 1+\Gf\Gj \\ 1+\Gf\Gj \epm,
\eeq{5.8}
that we will associate with the upper bulk modulus bounds, with
\beq  v_1=\frac{1-\Ga_2c_1}{\Ga_1\overline{c_1}}, \eeq{5.8a}
where the real constant $\Gf$ remains to be determined and $\Gj$ parameterizes the trajectory.
When $\Gj=0$ these are the fields in a rotation of the original crystal, and the term proportional to
$\Gj$ represents a stress jump, of the same form as $\BGs_1-\BGs_2$ in \eq{1.12}.

We next normalize and rotate the average fields, using the rotation
  \beq \BR=\bpm  e^{-i2\Gt}(1+e^{i\Gt}\Gj)/(1+e^{-i\Gt}\Gj) & 0 & 0 & 0 \\ 0 & e^{i2\Gt}(1+e^{-i\Gt}\Gj)/(1+e^{i\Gt}\Gj) & 0 & 0 \\ 0 & 0 & 1 & 0\\ 0 & 0 &0 & 1 \epm,
  \eeq{5.9}
to obtain
\beq \Bv'=\frac{\BR\BGs(\Gj)}{1+e^{-i\Gt}\Gj}=\bpm v'_1 \\ 1 \\ 1 \\ 1 \epm, \quad
\Bc'=\frac{\BR\BGs_0(\Gj)}{1+\Gf\Gj}=\bpm c_1' \\ \overline{c_1'} \\ 1 \\ 1 \epm,
\eeq{5.10}
with
\beq v_1'=\frac{(v_1+e^{-i3\Gt}\Gj)(1+e^{i\Gt}\Gj)}{(1+e^{-i\Gt}\Gj)^2},\quad
c_1' =\frac{(c_1+e^{-i2\Gt}\Gf\Gj)(1+e^{i\Gt}\Gj)}{(1+\Gf\Gj)(1+e^{-i\Gt}\Gj)}.
\eeq{5.11}
Substituting these in \eq{5.4} and using \eq{5.8a} and the relation $\Ga_1+\Ga_2=1$ gives
\beq \Gf=\frac{\overline{c_1}(1-\Ga_1\overline{c_1}e^{-i2\Gt}-\Ga_2c_1e^{i2\Gt})e^{-i3\Gt}}{1-\Ga_2c_1+\Ga_2\overline{c_1}e^{-i4\Gt}-\overline{c_1}e^{-i2\Gt}},
\eeq{5.12}
the reciprocal of which takes a simpler form:
\beq \frac{1}{\Gf}= \frac{e^{i\Gt}}{\overline{c_1}}+\frac{(e^{i\Gt}-e^{-i\Gt})(\overline{c_1} -e^{i\Gt})}{\Ga_1e^{-2i\Gt}(\overline{c_1} -e^{2i\Gt})+\Ga_2e^{2i\Gt}(c_1 -e^{-2i\Gt})}.
\eeq{5.12a}
Candidate values of $\Gt$ are determined by the requirement that $\Gf$ is real.

We now present an argument supplemented by numerical calculations of Christian Kern that convincingly indicates that point B is always attained when $\BC>0$.

 A trajectory $ c_1'(\Gj)$ is traced as $\Gj$ varies and this trajectory
 passes through the point $c_1$ at $\Gj=0$. In contrast to the trajectories \eq{3.19} and \eq{4.18} it is not generally the case that if $c_1'$ is on the trajectory
 then $\overline{c_1'}$ is also on the trajectory.
Consequently points where the trajectory self-intersects now do not typically correspond to orthotropic materials (having real values of $c_1'$). This is what
makes this case (point B) and the following case (point C) much more difficult to treat than the previous two cases (points A and D).

 For $\Gj$ large enough,
 $c_1'$ is close to $1$. To avoid redundancy we seek to
 identify each trajectory not by a point along it, but rather in terms
 of its limiting behavior near 1. 
 Let $\Ge=c_1-1$ and $\Gt$ be small parameters, such that the ratio
$\Ge/\Gt$ remains fixed as $\Ge\to 0$. Then from \eq{5.12a} we obtain, to first order in $\Gt$ and $\Ge$, that
\beq
\frac{1}{\Gf}\approx 1-(\overline{ z}-i)\Gt+
\frac{2i\Gt\overline{z}}{\Ga_1\overline{ z}+\Ga_2 z},\quad\text{where}\quad z=2i+(\Ge/\Gt).
\eeq{A2}
As $\Gf$ is real we get
\beqa
\frac{1}{\Gf} & = &  1-\Gt\Gd,\quad\text{with}\quad \Gd=z_R-\frac{4z_Rz_I}{z_R^2+(\Ga_1-\Ga_2)^2z_I^2}, \nonum
0 &= & z_I-1+\frac{2(z_R^2+(\Ga_1-\Ga_2)z_I^2)}{z_R^2+(\Ga_1-\Ga_2)^2z_I^2}. \eeqa{A3}
The latter can be solved for $z_R$ in terms of $z_I$ giving
\beq z_R^2=-\frac{[2+(z_I-1)(\Ga_1-\Ga_2)](\Ga_1-\Ga_2)z_I^2}{z_I+1}.
\eeq{A4}
Setting $\Gj=-1-t\Gt$, the terms entering the expression for $c_1'$ in 
\eq{5.11} are
\beqa c_1+e^{-i2\Gt}\Gf\Gj\approx \Gt(z-\Gd-t),\quad 1+e^{i\Gt}\Gj&\approx& -\Gt(i+t),\nonum
1+\Gf\Gj\approx -\Gt(t+\Gd),\quad 1+e^{-i\Gt}\Gj&\approx& -\Gt(-i+t),
\eeqa{A5}
to first order in $\Gt$. Substitution in \eq{5.11} gives an alternative
parameterization of each trajectory:
\beq  c_1' =\left(\frac{z}{\Gd+t}-1\right)\frac{1-it}{1+it}.
\eeq{A6}
A single trajectory is traced as $t$ varies, keeping $z$ fixed, and different values of $z_I$ generate different trajectories, with
$z_R$ being given by \eq{A4} ($z_I$ must be such that the right hand side of \eq{A4} is non-negative).
To see whether the trajectory $c_1'(t)$ loops around the origin as $t$ increases from $-\infty$ to $+\infty$ we look at
\beq \arg[c_1'(t)]=\arg[z/(\Gd+t)-1]+2\arg[1-it]
\eeq{A6a}
Now $\arg[z/(\Gd+t)-1]$ goes clockwise (or anticlockwise) through an angle of $\pi$ as $t$ increases from  $-\infty$ to $+\infty$ if $z_I$ is negative (positive).
Also $2\arg[1-it]$ goes clockwise through an angle of $2\pi$ as $t$ increases from  $-\infty$ to $+\infty$. So $c_1'(t)$ loops around the origin if
and only if these motions are both clockwise, i.e. if and only if $z_I<0$.

Additionally, we are only interested in those portions of the trajectory where $|c_1'(t)|<1$ since this necessarily holds if the corresponding elasticity
tensor $\BC'$ is positive definite. Thus we only need to consider $t$ in the interval $(|z|^2/(2z_R)-\Gd, \infty)$ if $z_R>0$, or in the interval
$(-\infty, |z|^2/(2z_R)-\Gd)$ if $z_R<0$.

%$0$ to $+\infty$ if $\Imag(e^{i2\Gt}c_1)$ is positive (or, respectively, negative), while $\arg[(e^{-i\Gt}\Gf-1)\Gl+1]$ goes anticlockwise (or clockwise) from
%$0$ through an angle of magnitude less than $\pi$ to $\arg[e^{-i\Gt}\Gf-1]$ as $\Gl$ goes from $0$ to $+\infty$ if $\Gf\sin(\Gt)$ is negative (or, respectively, positive). Thus $c_1'(\Gl)$ will loop around the origin as $\Gl$ varies from $0$ to $+\infty$ if the  total angle swept out is greater than or equal to $2\pi$, i.e. if
%\beq \arg[e^{i2\Gt}c_1-1]+2\arg[e^{-i\Gt}\Gf-1]\geq 2\pi\quad\text{or}\quad \arg[e^{i2\Gt}c_1-1]+2\arg[e^{-i\Gt}\Gf-1]\leq -2\pi.
%\eeq{5.15}

The points where the trajectory self intersects can be found by introducing $w=1/z=w_R+iw_I$
and taking $s=1/(\Gd +t)$ to parameterize the trajectory. The expression for $c_1'$ becomes
\beq c_1'=\frac{(w-s)[1+(i-\Gd)s]}{w[1-(i+\Gd)s]}. \eeq{A7}
Self intersections occur at $s=w_R\pm\Gg$ where $\Gg$ needs to be found from
\beq \frac{(iw_I-\Gg)[1+(i-\Gd)(w_R+\Gg)]}{[1-(i+\Gd)(w_R+\Gg)]}
=\frac{(iw_I+\Gg)[1+(i-\Gd)(w_R-\Gg)]}{[1-(i+\Gd)(w_R-\Gg)]}.
\eeq{A8}
Multiplying this by the product of the denominators we see that the imaginary parts cancel,
and from the real parts we get
\beq 1-2\Gd w_R+(\Gd^2+1)(w_R^2-\Gg^2)+2w_I=0, \eeq{A9}
giving
\beq \Gg=\pm\sqrt{w_R^2+\frac{1-2\Gd w_R+2w_I}{\Gd^2+1}}. \eeq{A10}

Numerical simulations of Christian Kern, which convincingly indicate that point B is always attained when $\BC>0$, are presented in \fig{5} and \fig{6}.
For various positive values of $\Ga_1$, they show the trajectories in the unit disk $|c_1|<1$, excluding the part of each trajectory
that loops around the origin. Each trajectory tail on one side of the loop is colored in red, while the trajectory tail on the other side
of the loop are colored in blue so that one can clearly see the point where the loop begins and ends. Only the trajectories tails in the lower
half of the  $c_1$-plane are shown: mirroring them about the real axis gives the trajectories tails in the upper half plane. The trajectory
tails, together with their mirrors, apparently fill the whole unit disk $|c_1|<1$, no matter what the value of $\Ga_1>0$. This strongly suggests that 
for any $\BC>0$ there exists a trajectory tail through $c_1$ and an associated point $c_1'$ on the trajectory where the
trajectory self intersects.

The one remaining question is whether the point $c_1'$ where the trajectory self intersects corresponds to a positive definite
tensor $\BC'$ when $c_1$ corresponds to a positive definite tensor $\BC$? One can argue as follows. First note that given  $c_1'$
at the self intersection point then assuming we can find one tensor $\BC'>0$ on the  that is associated with $c_1'$ then the
other elasticity tensor that is associated with $c_1'$ on the  can be taken to be a rotation of $\BC'$ by $2\Gt$.

Thus, under this assumption, associated with the trajectory tail and loop is
a hierarchical laminate polycrystal of the original crystal with a structure similar to that in \fig{1}.
The geometry of the hierarchical laminate is implied by our solution which gives hierarchical stress and strain fields in this medium. This is
a particular case of a microgeometry corresponding to the differential scheme and the tensor $\BC'$ can be found as the
solution to a ``self-consistent'' equation: $\BC'$ must be such that when this tensor is appropriately laminated with
the original crystal with tensor $\BC_0$, as dictated by the parameters associated with the trajectory,
the resulting laminate has an effective tensor which is a rotation of $\BC'$. The self-consistent equation can be obtained
from the lamination formula of Francfort and Murat \cite{Francfort:1986:HOB} or from one of the other lamination formulas described in chapter 9 of
\cite{Milton:2002:TOC}.
The realizability of  the differential scheme, proved for the original Bruggeman's differential scheme in \cite{Milton:1985:SEM} and in full generality in
\cite{Avellaneda:1987:IHD}, shows this ``self-consistent'' equation has a unique solution for $\BC'$ such that  $\BC'>0$ when $\BC_0>0$. 
There could be other solution branches where $\BC'$ is not positive definite, but we pick the one where  $\BC'>0$

%Now  $\BC$ depends on $c_1$ in some analytic way, and
%the corresponding $\BC'$ for each solution branch also will have an analytic dependence on $c_1$. One of these solution
%branches is the $\BC'(c_1')$ corresponding to the self-intersection point of the trajectory.
%Is it the unique solution branch such that $\BC'>0$ when $\BC>0$? When $c_1$ is very small, then $\BC$ will be positive definite as $c_1=0$ 
%corresponds to an isotropic $\BC$ with positive bulk and shear moduli inside the box $ABCD$. 
%The numerics show that then $c_1'$ is also very small, implying that $\BC'>0$.
%So it is the correct solution branch when  $c_1$ is very small, and by analytic continuation it must be the  correct solution branch
%for larger values of $c_1$ corresponding to crystals with $\BC>0$.

Having obtained a trajectory of optimal polycrystals including a loop around the origin associated with positive definite elastic tensors,
we can identify two optimal polycrystals with orthotropic symmetry corresponding to the two points where the loop intersects the real axis.
These orthotropic tensors typically have a different volume fraction in the first stage of the construction process, as illustrated in \fig{1}. Then, as a final step, we use the construction scheme in \cite{Avellaneda:1996:CCP} to obtain an optimal elastically isotropic material that corresponds
to the point $B$. 

\begin{figure}[!ht]
	\centering
	\begin{subfigure}{.50\linewidth}
		\includegraphics[width=\textwidth]{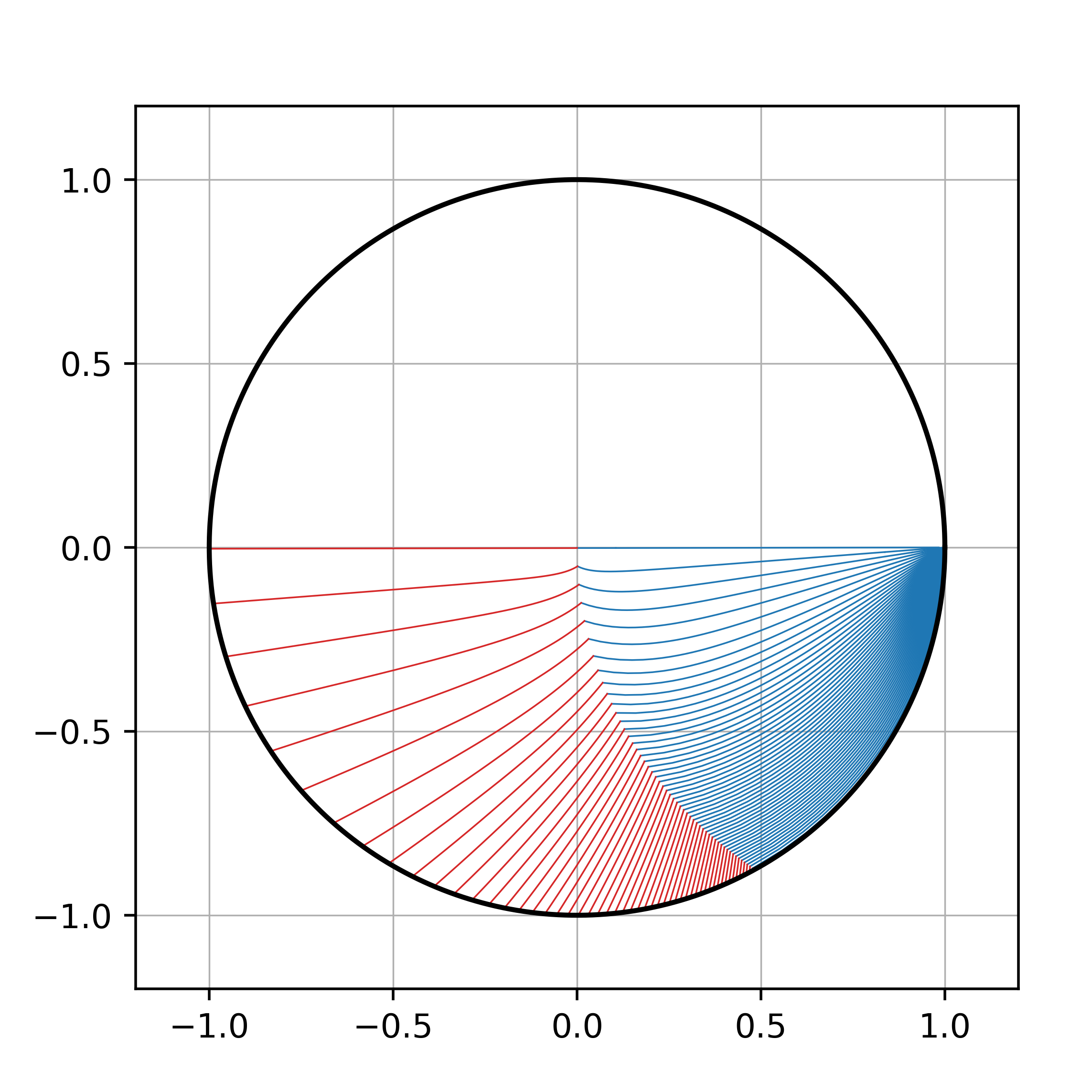}
		\caption{$\Ga_1=0.1$}
	\end{subfigure}
	\hskip-2em
	\begin{subfigure}{.50\linewidth}
		\includegraphics[width=\textwidth]{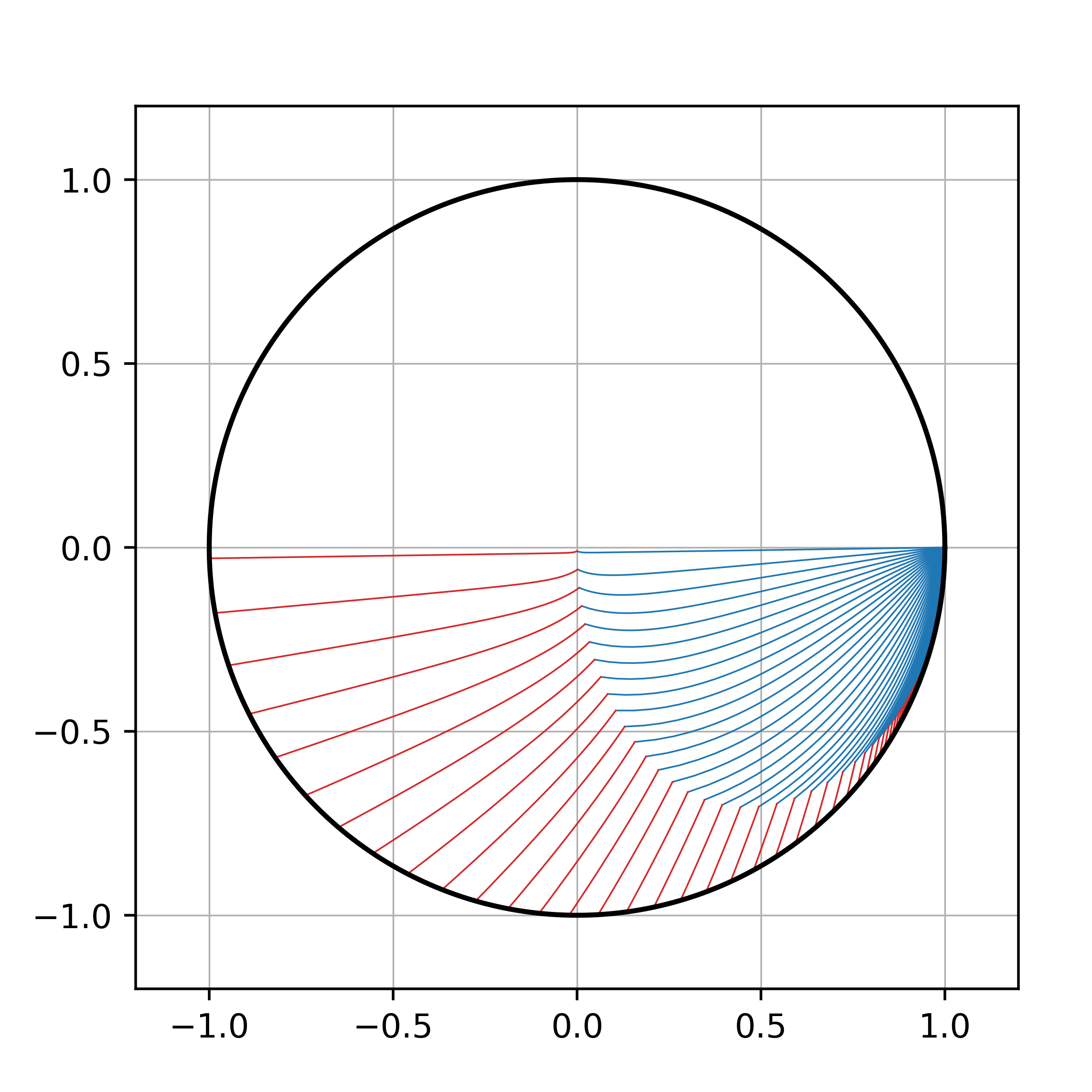}
		\caption{$\Ga_1=0.49$}
	\end{subfigure}
	\begin{subfigure}{.50\linewidth}
		\includegraphics[width=\textwidth]{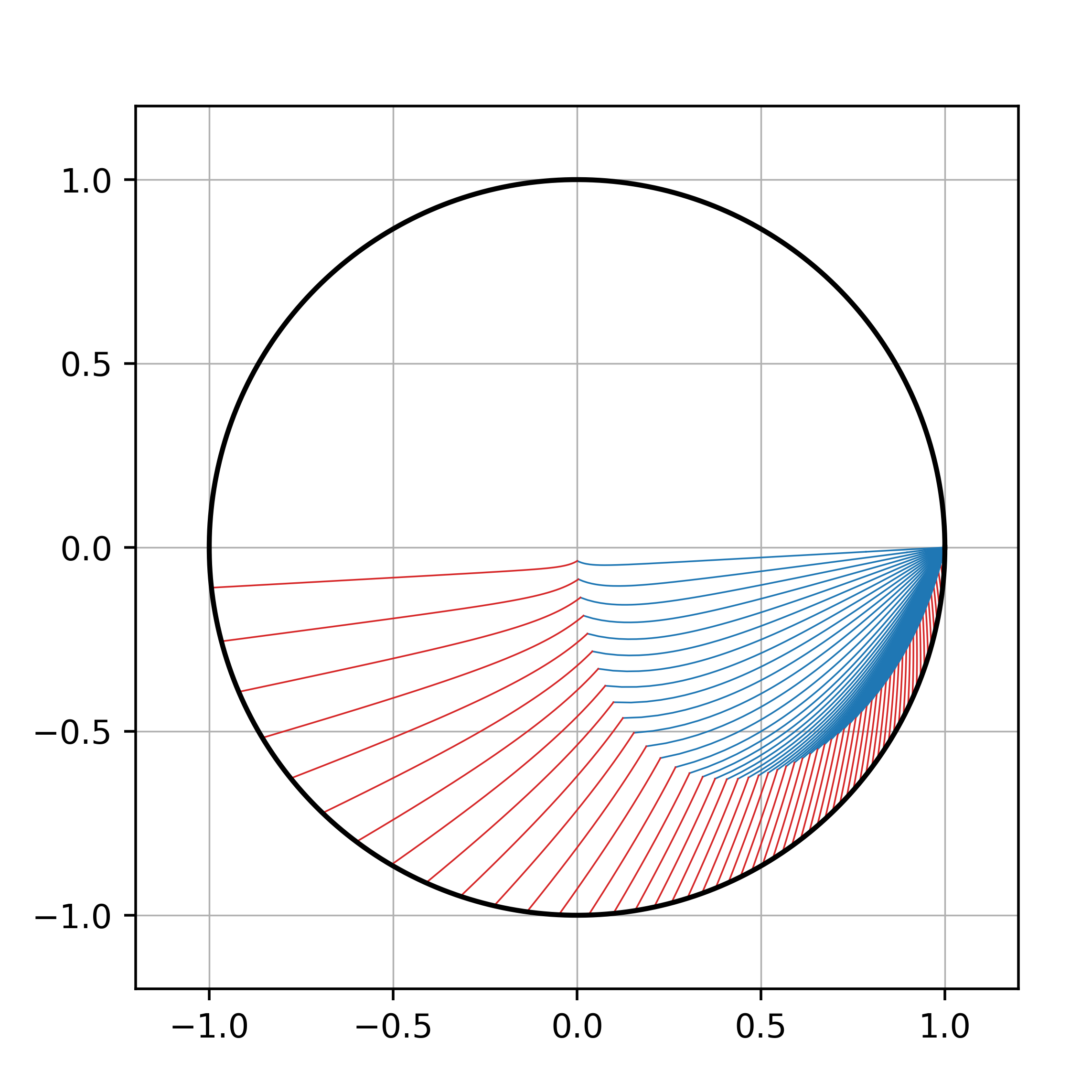}
		\caption{$\Ga_1=0.75$}
              \end{subfigure}
              	\hskip-2em
              	\begin{subfigure}{.50\linewidth}
		\includegraphics[width=\textwidth]{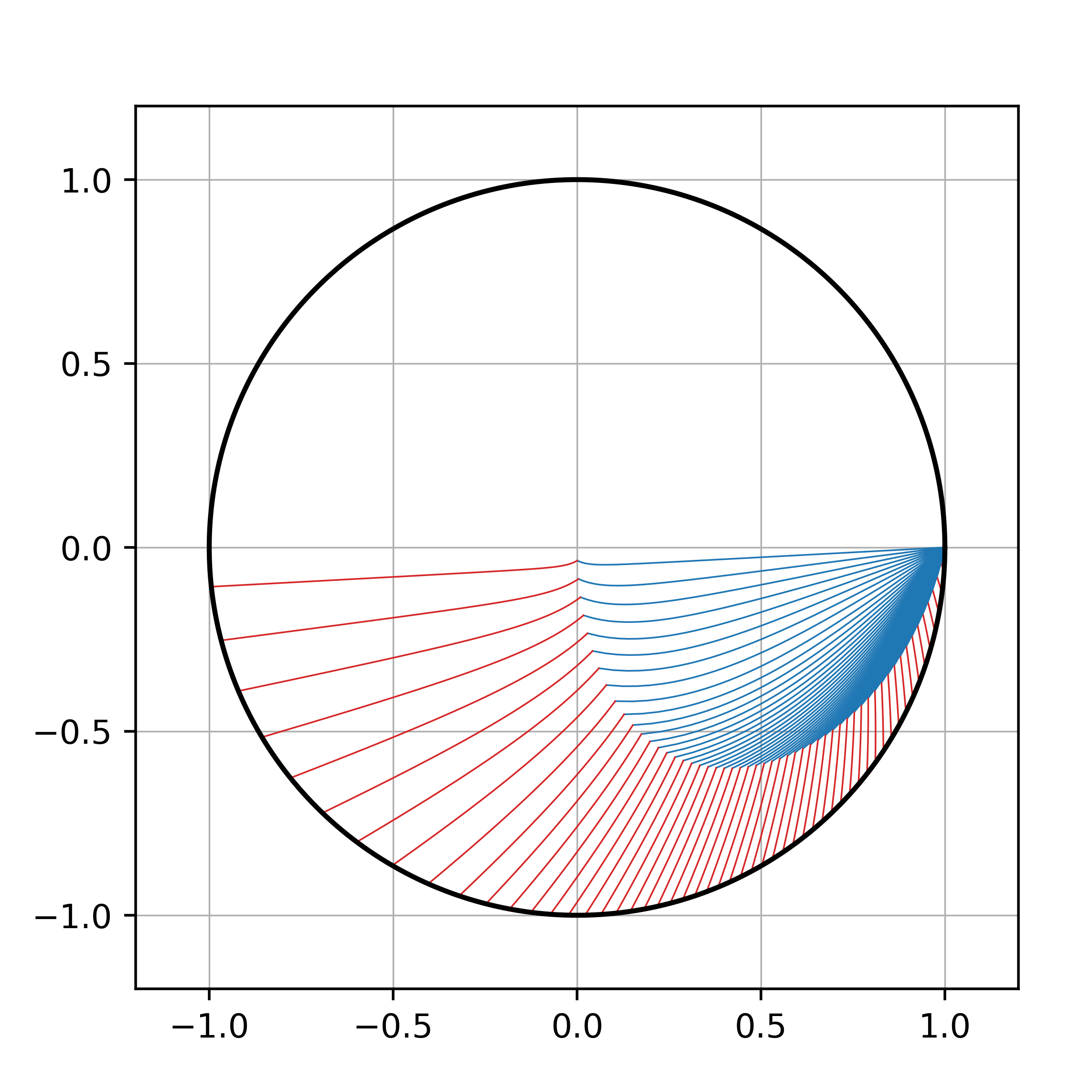}
		\caption{$\Ga_1=0.99$}
              \end{subfigure}
              \caption{Trajectory plots within the unit disk of the $c_1$-plane, excluding the loop portions, for some positive values of $\Ga_1<1$.
                Figure courtesy of Christian Kern.}
              \labfig{secondp}
              \labfig{5}
      \end{figure}

  \begin{figure}[!ht]
	\centering            
                \begin{subfigure}{.50\linewidth}
                  \includegraphics[width=\textwidth]{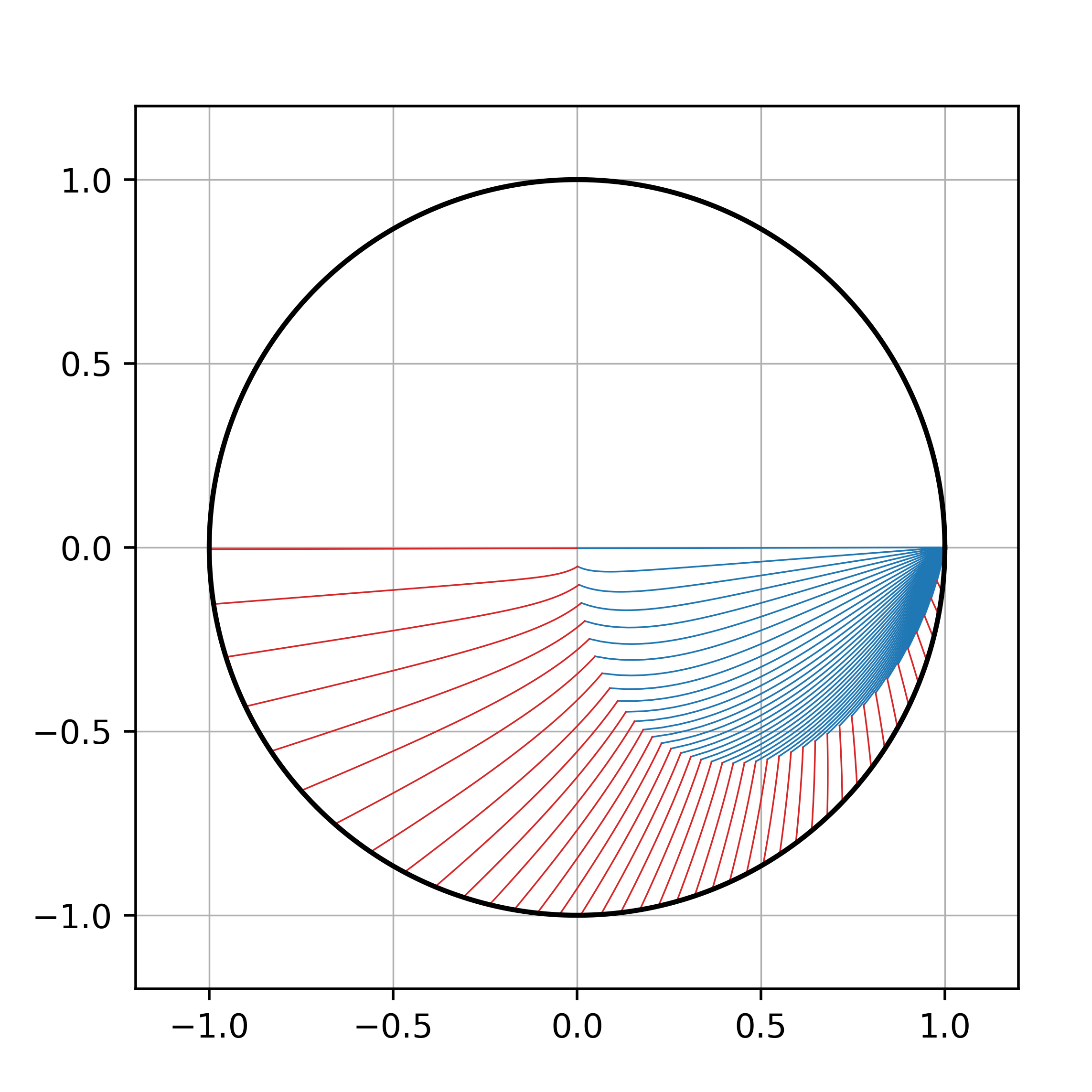}
		\caption{$\Ga_1=1.25$}
              \end{subfigure}
              	\hskip-2em
                \begin{subfigure}{.50\linewidth}
                  \includegraphics[width=\textwidth]{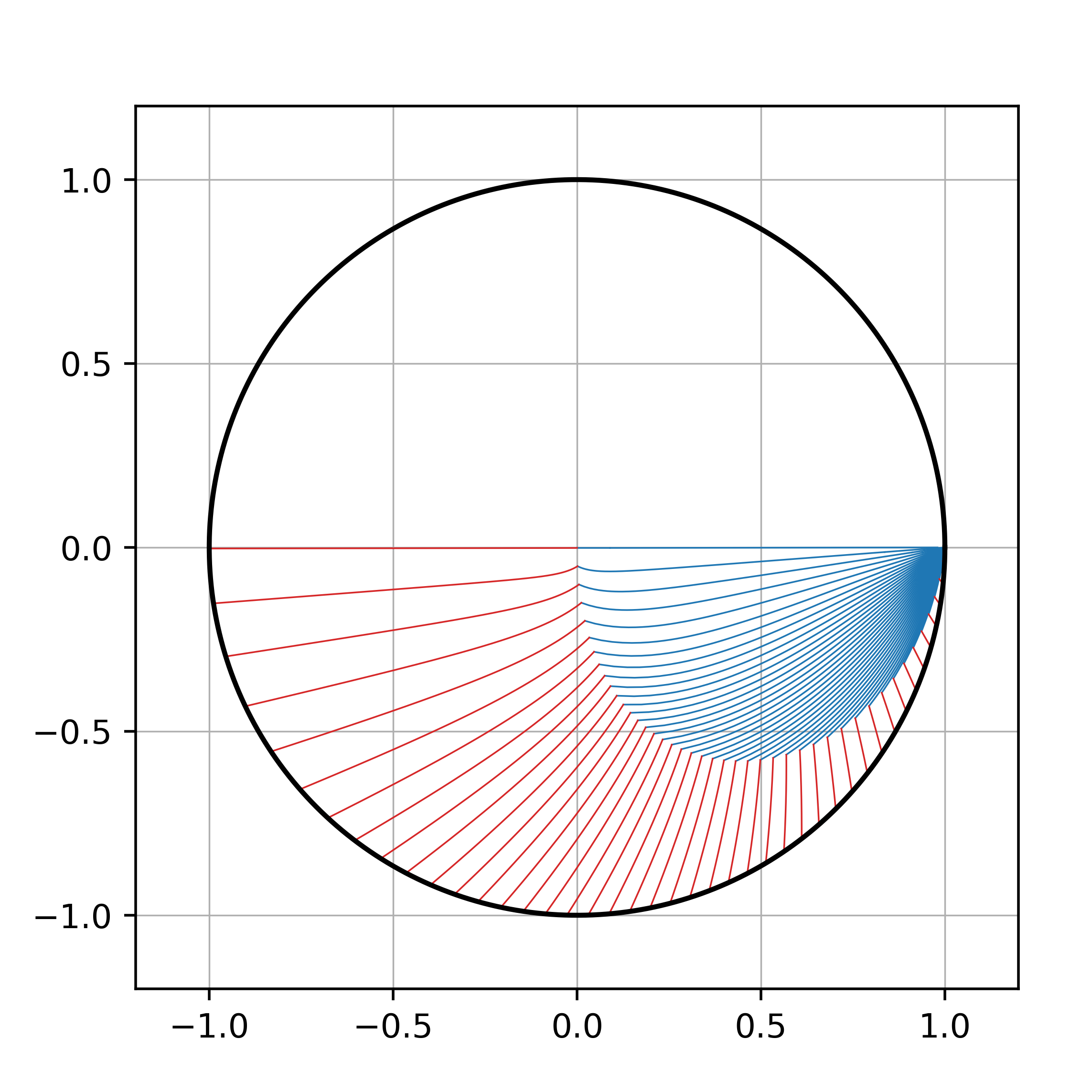}
                  \caption{$\Ga_1=1.51$}
              \end{subfigure}
                \begin{subfigure}{.50\linewidth}
                  \includegraphics[width=\textwidth]{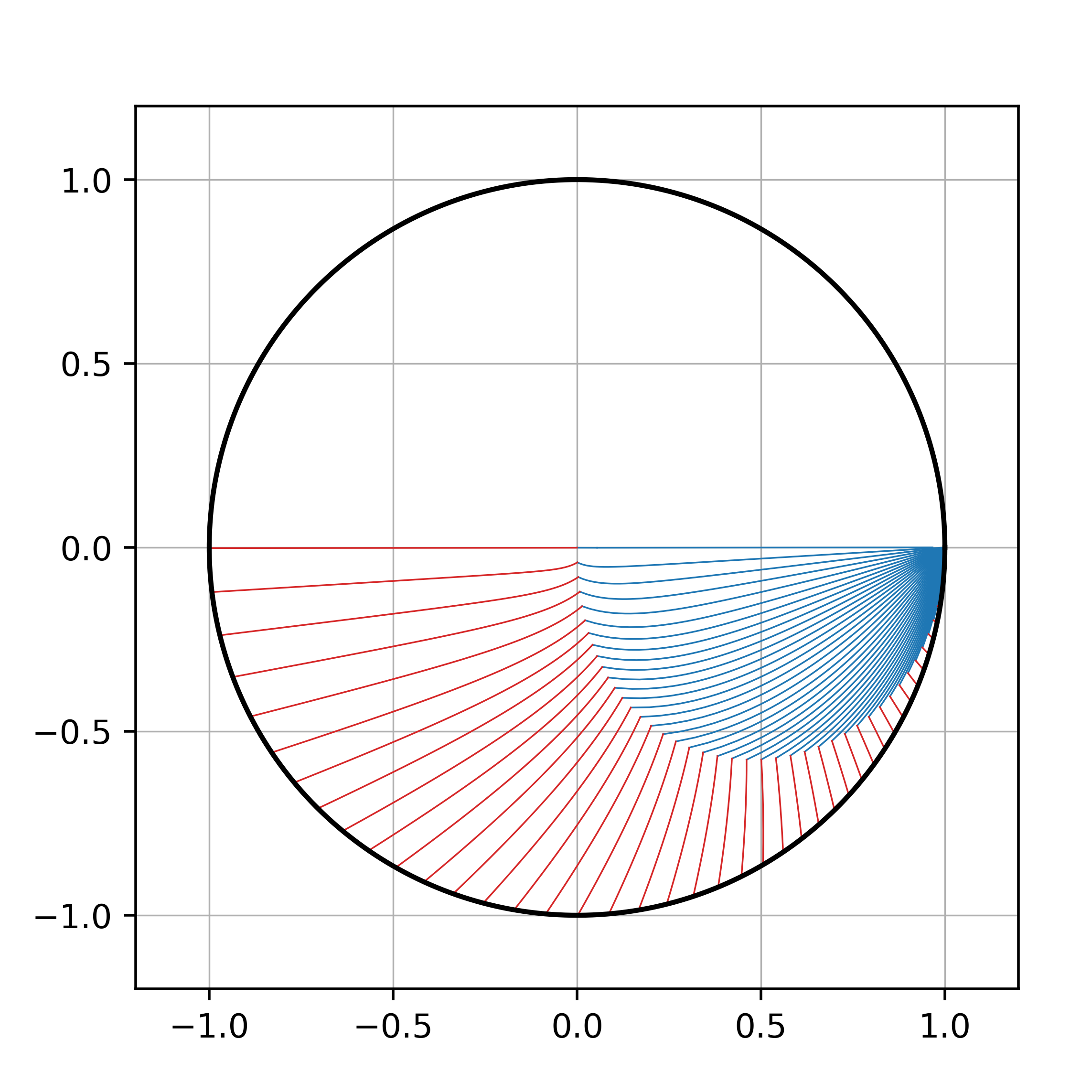}
		\caption{$\Ga_1=2$}
              \end{subfigure}
              	\hskip-2em
                             \begin{subfigure}{.50\linewidth}
                  \includegraphics[width=\textwidth]{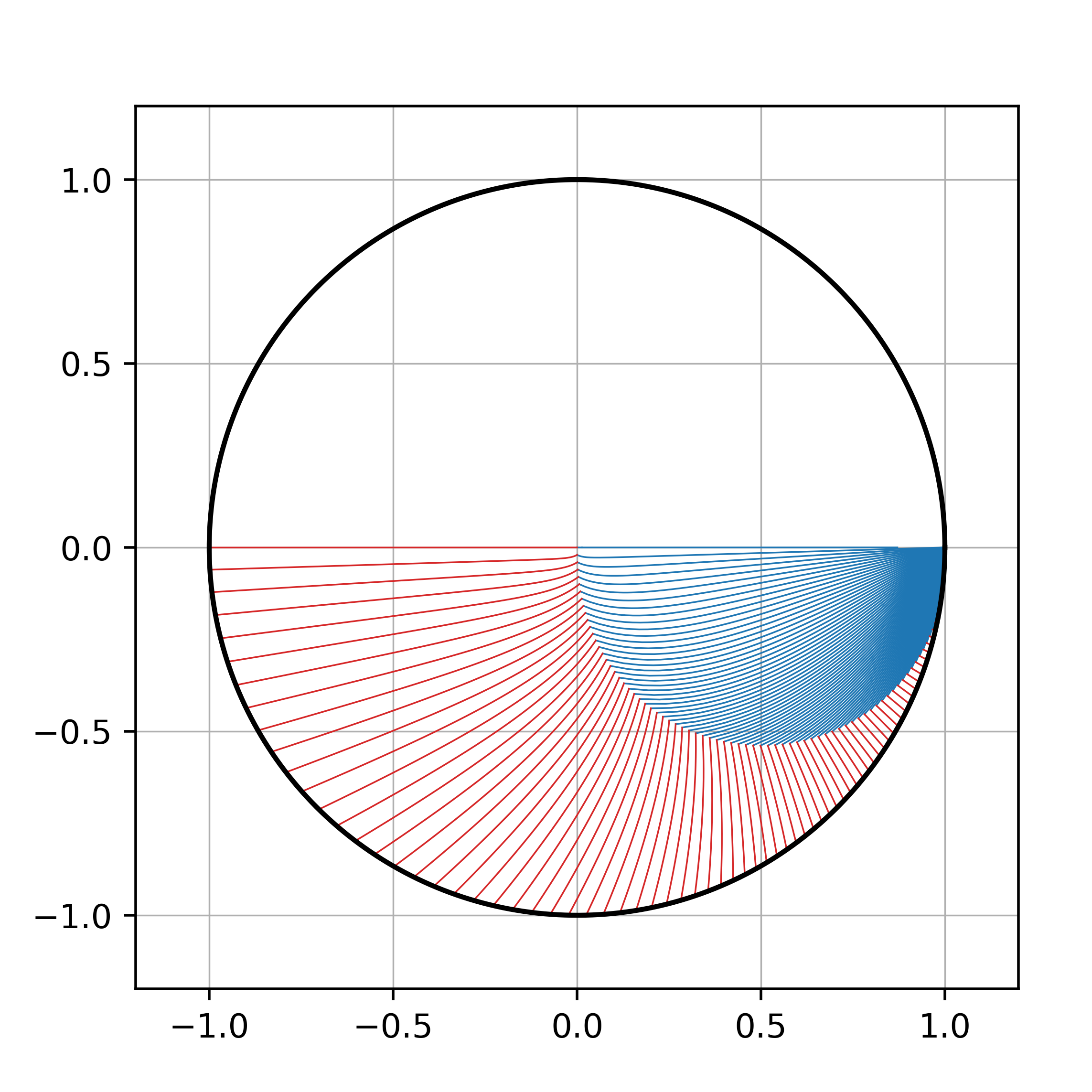}
		\caption{$\Ga_1=6$}
              \end{subfigure}
              \caption{Trajectory plots within the unit disk of the $c_1$-plane, excluding the loop portions,
                for some values of $\Ga_1>1$, corresponding to negative values of $\Ga_2=1-\Ga_1$. Figure courtesy of Christian Kern.}
	\labfig{6}
      \end{figure}

\newpage
 %%%%%%%%%%%%%%%%%%%%%%%%%%%%%%%%%%%%%%%%%%%%%%%%%%%%%%%%%%%%%%%%%%%%%%%%
\section{An algorithm for producing microstructures that simultaneously attain the upper bulk modulus and lower shear modulus bounds (point C)}
\setcounter{equation}{0}
%%%%%%%%%%%%%%%%%%%%%%%%%%%%%%%%%%%%%%%%%%%%%%%%%%%%%%%%%%%%%%%%%%%%%%%%%%
We look for geometries having a possibly anisotropic effective elasticity tensor $\BC'$ such that for some $\Bw'\ne 0$ and $\Bc'\ne 0$,
each representing symmetric matrices,
\beq \BC'-\BT\geq 0,\quad \BS'-\BT_0\geq 0,\quad\BC'\Bw'=\BT\Bw', \quad\Bc'=\BC'\BT_0\Bc'. \eeq{6.1}
Ultimately $\BC'$ will represent an effective
tensor of a polycrystal obtained from our starting material with elasticity tensor $\BC_0$. As $\BT\Bw'=\BC'\Bw'$ is a symmetric matrix, \eq{4.2} must hold.
We deduce that
\beq \overline{\Bc'}:\Bw'=\overline{\BC'\BT_0\Bc'}:\Bw'=\overline{\Bc'}:\BT_0\BT\Bw',
\eeq{6.2}
or equivalently, by normalization and by rotating and redefining $\BC'$ as necessary so that
\beq |w_1'|\geq t_2/t_1, \quad w_2'=-1, \quad w_3'=t_2/t_1, \quad w_4=1\quad \text{and} \quad c_2'=\overline{c_1'},\quad c_3'=c_4'=1, \eeq{6.3}
\eq{6.2} reduces to
\beq \Gb_1(t_1/t_2)\overline{c'_1}w'_1-\Gb_2c'_1=-1, \eeq{6.4}
where
\beq \Gb_1=\frac{t_2(1+t_0t_1)}{t_1+t_2 +2t_0t_1t_2},\quad\frac{t_1(1+t_0t_2)}{t_1+t_2 +2t_0t_1t_2}=1-\Gb_1 . \eeq{6.5}
Note that the last inequality in \eq{2.13a}, with $1/\Gk_*=2t_0$ (so that the bulk modulus is attained), implies that $t_1+t_2 +2t_0t_1t_2\geq 0$.

We can think of any material attaining the bounds as being parameterized by the complex number $c_1'$. In terms of it \eq{6.4} implies
\beq w_1'=t_2(\Gb_2c_1'-1)/(t_1\Gb_1\overline{c_1'}), \eeq{6.5a}
and the constraint that $|w_1'|\geq |w_3'|=|t_2/t_1|$ holds if and only if
\beq |\Gb_2c_1'-1|^2\geq \Gb_1^2|c_1'|^2.
\eeq{6.6}
Like \eq{5.6}, this is automatically satisfied if $|c_1'|\leq 1$.

Now consider the displacement gradient field trajectory
\beq \BE(\Gj)=\bpm e^{i2\Gt}w_1-t_2e^{-i\Gt}\Gj/t_1 \\  -e^{-i2\Gt}-e^{-i\Gt}\Gj  \\ t_2(1+e^{-i\Gt}\Gj)/t_1 \\ 1+e^{-i\Gt}\Gj \epm,
\eeq{6.7}
that we will associate with the upper shear modulus bounds, and the stress field trajectory
\beq \BGs_0(\Gj)=\bpm e^{i2\Gt}c_1+\Gf\Gj \\  e^{-i2\Gt}\overline{c_1}+\Gf\Gj  \\ 1+\Gf\Gj \\ 1+\Gf\Gj \epm,
\eeq{6.8}
that we will associate with the upper bulk modulus bounds, with
\beq  w_1=\frac{t_2(\Gb_2c_1-1)}{t_1\Gb_1\overline{c_1}}, \eeq{6.8a}
where the real constant $\Gf$ remains to be determined and $\Gj$ parameterizes the trajectory.
When $\Gj=0$ these are the fields in a rotation of the original crystal, and the term proportional to
$\Gj$ in \eq{6.7} represents a displacement gradient field jump, of the same form as $\BE_1-\BE_2$ in \eq{1.12}, while
in \eq{6.8} it represents a stress field jump, of the same form as $\BGs_1-\BGs_2$ in \eq{1.12}.

We next normalize and rotate the average fields, using the rotation
  \beq \BR=\bpm  e^{-i2\Gt}(1+e^{i\Gt}\Gj)/(1+e^{-i\Gt}\Gj) & 0 & 0 & 0 \\ 0 & e^{i2\Gt}(1+e^{-i\Gt}\Gj)/(1+e^{i\Gt}\Gj) & 0 & 0 \\ 0 & 0 & 1 & 0\\ 0 & 0 &0 & 1 \epm,
  \eeq{6.9}
to obtain
\beq \Bw'=\frac{\BR\BGs(\Gj)}{1+e^{-i\Gt}\Gj}=\bpm w'_1 \\ -1 \\ 1 \\ 1 \epm, \quad
\Bc'=\frac{\BR\BGs_0(\Gj)}{1+\Gf\Gj}=\bpm c_1' \\ \overline{c_1'} \\ 1 \\ 1 \epm,
\eeq{6.10}
with
\beq w_1'=\frac{(w_1-t_2e^{-i3\Gt}\Gj/t_1)(1+e^{i\Gt}\Gj)}{(1+e^{-i\Gt}\Gj)^2},\quad
c_1' =\frac{(c_1+e^{-i2\Gt}\Gf\Gj)(1+e^{i\Gt}\Gj)}{(1+\Gf\Gj)(1+e^{-i\Gt}\Gj)}.
\eeq{6.11}
Substituting these in \eq{6.4} and using \eq{6.8a} and the relation $\Gb_1+\Gb_2=1$ gives
 \beq \Gf=\frac{\overline{c_1}(1-\Gb_1\overline{c_1}e^{-i2\Gt}t_1/t_2-\Gb_2c_1e^{i2\Gt})e^{-i3\Gt}}{1-\Gb_2c_1+\Gb_2\overline{c_1}e^{-i4\Gt}-\overline{c_1}e^{-i2\Gt}}.
\eeq{6.12}
Candidate values of $\Gt$ are determined by the requirement that $\Gf$ is real. Note that the expression \eq{6.11} for $c_1'$ is exactly the same as in \eq{5.11}
and consequently the condition for the trajectory $c_1'(\Gj)$ to loop around the origin is the same as that given in the previous section. The barrier
to numerically testing for realizability of point $C$ for all $\BC>0$ is that we need to analyze the trajectories in the $c_1$-plane as not just one but two parameters are varied
($\Gb_1=1-\Gb_2$ and $t_1/t_2$) since these both enter \eq{6.12}. However, given a specific $\BC>0$ the algorithm can be easily implemented and if one finds
that there is a trajectory tail through the corresponding $c_1$ with the trajectory looping around the origin and self-intersecting at some $c_1'$, corresponding
to a tensor $\BC'>0$, then realizability of point $C$ will be established for that  $\BC$. 

\section*{Acknowledgements}
The author is grateful to the National Science Foundation for support through the Research Grant DMS-1814854.
Christian Kern is thanked for the numerical
simulations presented in Figures 5 and 6. Additionally, the author is deeply grateful to one referee who carefully examined the paper and noted many
points requiring correction. In particular they noticed a significant gap in the original argument for attainability of points A and D, which is now corrected.

%%\bibliographystyle{../siamplain}
%\bibliography{/u/ma/milton/tcbook,/u/ma/milton/newref}
%\bibliography{/home/milton/tcbook,/home/milton/newref}

\ifx \bblindex \undefined \def \bblindex #1{} \fi\ifx \bbljournal \undefined
  \def \bbljournal #1{{\em #1}\index{#1@{\em #1}}} \fi\ifx \bblnumber
  \undefined \def \bblnumber #1{{\bf #1}} \fi\ifx \bblvolume \undefined \def
  \bblvolume #1{{\bf #1}} \fi\ifx \noopsort \undefined \def \noopsort #1{}
  \fi\ifx \bblindex \undefined \def \bblindex #1{} \fi\ifx \bbljournal
  \undefined \def \bbljournal #1{{\em #1}\index{#1@{\em #1}}} \fi\ifx
  \bblnumber \undefined \def \bblnumber #1{{\bf #1}} \fi\ifx \bblvolume
  \undefined \def \bblvolume #1{{\bf #1}} \fi\ifx \noopsort \undefined \def
  \noopsort #1{} \fi

\end{document}